\newif\ifpdf\ifx\pdfoutput\undefined\pdffalse%
                \documentclass{article}\usepackage{graphicx}%
        \else\pdfoutput=1\pdftrue%
                \documentclass[pdftex]{article}%
                \usepackage[pdftex]{hyperref}%
                \usepackage[pdftex]{graphicx}%
        \fi

\usepackage{amsmath}
\usepackage{amsfonts}
\usepackage{amssymb}
\usepackage[normalem]{ulem}

\newtheorem{theorem}{Theorem}[section]

\newtheorem{corollary}[theorem]{Corollary}

\newtheorem{lemma}[theorem]{Lemma}

\newtheorem{remark}[theorem]{Remark}
\newenvironment{proof}[1][Proof]{\textbf{#1.} }%
               {\ \hfill\rule{0.5em}{0.5em}\medskip}

\newcommand{\qed}{\hfill\rule{0.5em}{0.5em}\medskip}%
\numberwithin{equation}{section} 

\newif\ifShowComments
\ShowCommentstrue
\def\strutdepth{\dp\strutbox}
\def\druk#1{\strut\vadjust{\kern-\strutdepth
        {\vtop to \strutdepth{%
                \baselineskip\strutdepth\vss
                        \llap{\hbox{#1}\quad}\null}}}}

\def\bl{\bigl}%
\def\br{\bigr}%
\def\Bl{\Bigl}%
\def\Br{\Bigr}%
\def\cst{{\sf const}}%
\def\cstL{{\sf const}_L}%

\def\INp{\mathcal{I}_{N,+}}%
\newcommand{\ZNpl}[1]{Z_{N,+,\lambda}^{\,#1}}%
\newcommand{\Zpl}[2]{Z_{#1,+,\lambda}^{\,#2}}%
\newcommand{\Zpo}[2]{Z_{#1,+,0}^{\,#2}}%
\newcommand{\Pp}[2]{\mathrm{P}_{#1,+}^{\,#2}}%
\newcommand{\Ppo}[2]{\mathrm{P}_{#1,+,0}^{\,#2}}%
\newcommand{\Ppl}[2]{\mathrm{P}_{#1,+,\lambda}^{\,#2}}%
\newcommand{\PNpl}[1]{\mathrm{P}_{N,+,\lambda}^{\,#1}}%
\newcommand{\ENpl}[1]{\mathrm{P}_{N,+,\lambda}^{\,#1}}%
\newcommand{\epl}{\mathrm{E}_{+,\lambda}}%
\newcommand{\ppl}{\mathrm{P}_{+,\lambda}}%
\newcommand{\ZNpo}[1]{Z_{N,+,0}^{\,#1}}%
\newcommand{\PNpo}[1]{\mathrm{P}_{N,+,0}^{\,#1}}%
\newcommand{\Pmpo}[1]{\mathrm{P}_{\calI_m,+,0}^{\,#1}}%
\newcommand{\PDpo}[1]{\mathrm{P}_{\Delta,+,0}^{\,#1}}%
\newcommand{\ENpo}[1]{\mathrm{E}_{N,+,0}^{\,#1}}%
\newcommand{\Zl}[2]{Z_{#1,\lambda}^{\,#2}}%
\newcommand{\Pl}[2]{\mathrm{P}_{#1,\lambda}^{\,#2}}%
\def\calA{\mathcal{A}}%
\def\calB{\mathcal{B}}%
\def\calC{\mathcal{C}}%
\def\calD{\mathcal{D}}%
\def\veps{\varepsilon}%
\def\sig{\sigma}%
\def\neps{n_\varepsilon}%
\def\calI{\mathcal{I}}%
\def\calK{\mathcal{K}}%
\def\calJ{\mathcal{J}}%
\def\calE{\mathcal{E}}%
\def\Ad{\calA_\delta}%
\newcommand{\one}{\hbox{1\kern-.27em I}}%
\newcommand{\ind}[1]{\one_{#1}}
\def\Cov{{\sf Cov}}%
\def\every{{}^\forall}%
\def\some{{}^\exists}%
\def\BbbP{\mathbb{P}}%
\def\BbbR{\mathbb{R}}%
\def\BbbV{\mathbb{V}}%
\def\BbbX{\mathbb{X}}%
\def\BbbZ{\mathbb{Z}}%
\def\calN{\mathcal{N}}%
\def\Hgam{H_\gamma}%
\newcommand{\Hg}[1]{H_{\gamma_{#1}}}%
\def\Hone{H_1}%
\newcommand{\PPNpl}[1]{\BbbP_{N,+,\lambda}^{\,#1}}%
\newcommand{\PPNpo}[1]{\BbbP_{N,+,0}^{\,#1}}%
\newcommand{\PPmpo}[1]{\BbbP_{\calI_m,+,0}^{\,#1}}%
\def\nul{}
\def\punkt{\,\cdot\,}%
\def\ppunkt{\,\cdot\,,\,\cdot\,}%
\def\abx{a_x,b_x}%
\def\aby{a_y,b_y}%
\def\abxy{a_x,b_x,a_y,b_y}%
\def\abxym{a_x^m,b_x^m,a_y^m,b_y^m}%
\def\NN{\calN_{(0,N)}}%
\def\Nm{\calN_{\calI_m^o}}%
\def\noNm{\overline\calN_{\calI_m^o}}%
\def\Am{\calA_m}%
\def\Amup{\Am^{\nearrow}}%
\def\Amdn{\Am^{\searrow}}%
\def\Bm{\calB_m}%
\def\Bmup{\Bm^{\nearrow}}%
\def\Bmdn{\Bm^{\searrow}}%
\def\CmD{C^-_\Delta}%
\def\CpD{C^+_\Delta}%
\def\Htil{\widetilde H}%
\def\Deltil{\widetilde\Delta}%
\def\netil{\tilde n_\veps}%
\def\frakl{\mathfrak l}%
\def\frakr{\mathfrak r}%
\def\Delrho{\Delta_\rho}%
\def\hatD{\widehat D}%
\def\hatDom{\hatD^o_m}%
\def\rmP{\mathrm{P}}%
\def\comp{{\,c}}%
\def\sigtil{\tilde\sig}%

\def\bfE{{\mathbf E}}%
\def\bfP{{\mathbf P}}%

\title{Universality of Critical Behaviour in a Class of Recurrent
Random Walks}
\author{O.Hryniv\\Statistical Laboratory, DPMMS,\\ University of
Cambridge, Cambridge CB3~0WB, UK\\{\tt o.hryniv@statslab.cam.ac.uk }
\and 
Y.Velenik\\Laboratoire de Math\'ematiques Rapha\"el
Salem\\UMR-CNRS 6085, Universit\'e de Rouen\\
F-76821 Mont Saint Aignan\\{\tt Yvan.Velenik@univ-rouen.fr}}

\begin{document}
\maketitle
\begin{abstract}
Let $X_0=0$, $X_1$, $X_2$, \dots\, be an aperiodic random walk
generated by a 
sequence $\xi_1$, $\xi_2$, \dots\, of i.i.d.\ integer-valued random
variables with common distribution $p\,(\cdot)$ having zero mean and finite
variance. For an $N$-step trajectory $\BbbX=(X_0,X_1,\dots,X_N)$ and 
a monotone convex function $V:\BbbR^+\to\BbbR^+$ with $V(0)=0$, define
$\BbbV(\BbbX)=\sum_{j=1}^{N-1}V\bl(|X_j|\br)$. 
Further, let $\mathcal{I}_{N,+}^{a,b}$ be the set of all non-negative
paths $\BbbX$ compatible with the boundary conditions $X_0=a$,
$X_N=b$. 
We discuss asymptotic properties of $\BbbX\in\mathcal{I}_{N,+}^{a,b}$
w.r.t.\ the probability distribution 
\[
\Ppl{N}{a,b}(\BbbX)= \bl(\Zpl{N}{a,b}\br)^{-1}\,
\exp\Bl\{-\lambda\,\BbbV(\BbbX)\Br\}\; 
       \prod_{i=0}^{N-1} p\,(X_{i+1}-X_i)
\]
as $N\to\infty$ and $\lambda\to0$, $\Zpl{N}{a,b}$ being the
corresponding normalization.

If $V(\punkt)$ grows not faster than polynomially at infinity, define
$H(\lambda)$ to be the unique solution to the equation
\[
\lambda H^2\,V(H) =1\,.
\]
Our main result reads that as $\lambda\to0$, the  
typical height of $X_{[\alpha N]}$ scales as $H(\lambda)$ and the
correlations along $\BbbX$ decay exponentially on the scale
$H(\lambda)^2$.
Using a suitable blocking argument, we show that the distribution
tails of the rescaled height decay exponentially with critical
exponent~$3/2$.  
In the particular case of linear potential $V(\punkt)$, the
characteristic length $H(\lambda)$ is proportional to
$\lambda^{-1/3}$ as $\lambda\to0$.
\end{abstract}

\section{Introduction}
In this work, we are interested in the path-wise behaviour of a general class of random walks on the integers, whose path measure is submitted to a special form of exponential perturbation, the physical motivation of which is discussed at the end of this section. More precisely, to each $i\in\BbbZ$, we associate an integer\footnote
{  
Although in the sequel we'll discuss mainly integer-valued
one-dimensional random walks, analogous results for real-valued walks
can be obtained in a similar way. 
}
non-negative value
$X_i$ and for any integer interval 
\begin{equation}\label{eq:DeltaLR.def}
\Delta_{l,r}=(l,r)\equiv\bl\{l+1,l+2,\dots,r-1\br\}\subset\BbbZ
\end{equation}
we denote by $\mathcal{I}_{\Delta_{l,r},+}$ the
set of all such trajectories in $\Delta_{l,r}$: 
\[
\mathcal{I}_{\Delta_{l,r},+}=\mathcal{I}_{(l,r),+}
=\Bl\{\BbbX=(X_i)_{l<i<r}:X_i\ge0\Br\}. 
\]
Let $V:\BbbR_+\to\BbbR_+$ be a convex increasing continuous function
with $V(0)=0$ and a bounded growth at infinity:
\begin{quote}\sl
There exists 
$f:\BbbR^+\to\BbbR^+$ such that for any
$\alpha>0$ we have
\end{quote}
\begin{equation}\label{eq:Vregularity}
\limsup_{x\to\infty}\frac{V(\alpha x)}{V(x)}\le f(\alpha)<\infty.
\end{equation}
This property holds clearly for any (convex) polynomial function.

The probability of a trajectory $\BbbX\in\mathcal{I}_{(l,r),+}$ is defined
then via 
\begin{equation}\label{eq:PDpl.Def}
\Ppl{(l,r)}{a,b}(\BbbX) = \bl(\Zpl{(l,r)}{a,b}\br)^{-1}\,
\exp\Bl\{-\lambda\, \sum_{i=l+1}^{r-1} V\bl(X_i\br)\Br\}\; 
       \prod_{i=l}^{r-1} p(X_{i+1}-X_i)\,,
\end{equation}
where the boundary conditions are given by $X_l=a$ and $X_r=b$,
the parameter $\lambda$ is some strictly positive real number and 
$p(\punkt)$ are the transition probabilities of a 1D integer-valued
random walk with zero mean and finite second moment.
We suppose that the random walk is {\it strictly aperiodic\/} in the
sense that its $n$-step transition probabilities $p^n(\punkt)$ possess
the following property:
\begin{quote}
{\sl there is $A>0$ such that }
\end{quote}
\begin{equation}\label{eq:aperiodic}
\min\Bl\{p^n(-1),p^n(0),p^n(1)\Br\}>0 \quad \text{\sl for all\/ }n\ge A.
\end{equation}
When the boundary conditions are chosen such that $a=b=0$, we omit
them from the notation.
We will denote by $\Pl{\Delta_{l,r}}{a,b}$ the analogous probability
measure without the positivity constraint, and by
$\Zl{\Delta_{l,r}}{a,b}$ the associated partition function.
The special case of $\Delta_{0,N}=(0,N)$ will be abbreviated to
\begin{equation}\label{eq:INp.Def}
\INp=\Bl\{\BbbX=(X_i)_{1\leq i\leq N-1}:X_i\ge0\Br\}
\end{equation}
and
\begin{equation}\label{eq:PNpl.Def}
\PNpl{a,b}(\BbbX) = \bl(\ZNpl{a,b}\br)^{-1}\,
\exp\Bl\{-\lambda\, \sum_{i=1}^{N-1} V\bl(X_i\br)\Br\}\; 
       \prod_{i=0}^{N-1} p(X_{i+1}-X_i)\,
\end{equation}
respectively.

\bigskip
When $\lambda=0$ it is expected that, after a suitable rescaling, the
law of the random walk converges to that of a Brownian excursion; 
at least, it is a direct corollary ot the results in \cite{eB76}
and \cite{wdK76} that, in the diffusive scaling, the large-$N$ limit
of the distribution $\PNpo{}(\punkt)$ is non-trivial.
In particular, the path
delocalizes. When $\lambda>0$, the behaviour changes drastically:
We'll prove below that the path remains localized, and that the
correlations between positions $X_i$ and $X_j$ of the random walk
decay exponentially with their separation $|i-j|$. Our main goal is to
investigate how delocalization occurs as $\lambda$ decreases to
$0$. The corresponding critical behaviour can be analyzed in quite
some details and, most interestingly, under very weak assumptions on
the original random walk. In this way, it is possible to probe its
degree of universality. 

\bigskip
The physical motivation for the path measure considered in this work
is the phenomenon of critical prewetting. Consider a vessel containing
the thermodynamically stable gaseous phase of some substance. When the
boundary of the vessel displays a sufficiently strong preference
towards the thermodynamically unstable liquid phase, there may be
creation of a microscopic film of liquid phase coating the walls. As
the system is brought closer and closer to liquid/gas phase
coexistence, the layer of unstable phase starts to grow. For systems
with short-range interactions, two kind of behaviours are possible:
either there is an infinite sequence of first-order (so-called
layering) phase transitions, at which the thickness increases by one
mono-layer, or the growth occurs continuously; this is the case of
critical prewetting, and it is typical in two-dimensional systems, as
those modelled in the present work. We are thus interested in
quantifying the growth as a function of the distance to phase
coexistence. A natural parameter is the difference between the free
energy densities of the stable and unstable phases. Choosing
$V(x)=|x|$, we see that the perturbation $\lambda \sum_i V(X_i)$ can
be interpreted as the total excess free energy associated to the
unstable layer, the parameter $\lambda$ playing the role of the excess
free energy density. 

The problem of critical prewetting in continuous effective interface
models in higher dimensions, as well as in the 2D Ising model, has
been considered in~\cite{Ve2003}. The latter results are however
restricted to a much smaller class of interactions, and take also a
weaker form than those we obtain here. Notice however that the
thickness of the layer of the unstable phase in the 2D Ising model also
grows with exponent $1/3$, showing (not surprisingly) that this model
is in the same universality class as those considered in the present
work; one can hope therefore that the finer estimates we obtain here
have similar counterparts in the 2D Ising
model. Results  about critical
prewetting for one-dimensional effective interface models have already
been obtained in~\cite{dbAerS86}; they are limited to the
particular case of real-valued random walks with $V(x)=|x|$ and
$p(x)\propto e^{-\beta|x|}$, which turn out to be exactly
solvable. They are thus able to extract more precise information than
those we present here, including prefactors. It is not clear however
to what extent these finer properties also have universal
significance.

\smallskip
{\bf Acknowledgements.} 
It is a pleasure to thank H.Spohn for his interest in this work.
O.H.~also thanks L.C.G.Rogers, S.Molchanov and G.BenArous for stimulating
discussions.

\subsection{Main results}

We recall that the potential function $V(\punkt)$ in
\eqref{eq:PNpl.Def} is assumed to be continuous and increasing from
$0$ to $\infty$ as $x$ varies from $0$ to $\infty$. Therefore, for any
$\gamma>0$ there is a unique solution $\Hgam>0$ to the equation
\begin{equation}\label{eq:Hgam.Def}
\lambda\, H^2\,V(2\gamma H)\mid_{H=\Hgam}=1.
\end{equation}
The scale $\Hgam=\Hgam(\lambda)$ will play an important role in our
future considerations. As a simple corollary of the definition
\eqref{eq:Hgam.Def}, note that $\gamma^{1/3}\Hgam$ is a non-increasing
function of $\gamma$; indeed, thanks to convexity and monotonicity
of $V(\punkt)$, for any $0<\gamma_1<\gamma_2$ and any $0<a<1$ we get:
\[
V(2\gamma_1\Hg1)\Hg1^2\equiv V(2\gamma_2\Hg2)\Hg2^2=\frac1\lambda\ge 
V\bl(2a^3\gamma_2\Hg2/a\br)\bl(\Hg2/a\br)^2;
\]
now take $a$ satisfying $a^3\gamma_2=\gamma_1$ and recall monotonicity
of the function $x\mapsto V(2\gamma_1x)x^2$ to infer $a\Hg1\ge\Hg2$.

Our first result says that the ``average height'' of the interface in
the limit of small $\lambda$ is of order $\Hone(\lambda)$.

\begin{theorem}\label{thm:L1Bounds}\sl
Let $H=\Hone(\lambda)$ be as defined in \eqref{eq:Hgam.Def}. 
There exist positive constants  $\delta_0$,
$\lambda_0$,  and $C_1$, $C_2$ such that the inequalities 
\begin{gather}
\PNpl{}\Bl(\sum_{i=1}^NX_i \ge \delta^{-1}HN\Br) \leq
\frac1{C_1}e^{-C_1\delta^{-1}H^{-2}N}\,,
\label{eq:UpperBound}
\\
\PNpl{}\Bl(\sum_{i=1}^NX_i \le \delta\,HN\Br) \leq
\frac1{C_2}e^{-C_2\delta^{-2}H^{-2}N}
\label{eq:LowerBound}
\end{gather}
hold uniformly in $0<\lambda<\lambda_0$, $0<\delta\le\delta_0$, and
$N\ge H^2$. 
\end{theorem}

Our next result describes the tails of the point-wise height
distribution. Although it is formulated for the height in the
middle of a typical interface, the result holds for all points
from $[AH^2,N-AH^2]$, with any fixed $A>0$, ie., lying sufficiently
deep in $[0,N]$; then, clearly, $c_i=c_i(A)\to0$, $i=1,2$, as $A\to0$.

\begin{theorem}\label{thm:PointwiseBounds}\sl
Let $H=\Hone(\lambda)$ be as defined in \eqref{eq:Hgam.Def}. There
exist positive constants $T_0$, $K$, $c_1$ and $c_2$ such
that for any $T\ge T_0$ and all $N\ge KH^2$ the inequalities
\begin{equation}\label{eq:PointwiseBounds}
\frac1{c_1}\exp\Bl\{-c_1T^{3/2}\Br\}\le
\PNpl{}\Bl(X_{[N/2]}>T\Hone\Br)\le 
\frac1{c_2}\exp\Bl\{-c_2T^{3/2}\Br\}
\end{equation}
hold for all $\lambda\in(0,\lambda_0]$, where
$\lambda_0=\lambda_0(T)>0$. 
In particular, these estimates are uniform on compact subsets of
$[0,\infty)$. 
\end{theorem}

\begin{remark}\sl
Tail estimates for small $\lambda$ uniform in $T$ large enough can be
obtained taking into account the tail behaviour of the original
random-walk. For example, a variant of the argument above, actually
even simpler, proves that for a Gaussian random walk,
i.e. $p(x)\propto e^{-c|x|^2}$, the tail behaviour proved above holds
with the same exponent for all $\lambda\leq \lambda_0$ uniformly in
large $T$. This behaviour is not universal however, as different
behaviour of $p(\,\cdot\,)$ can give rise to completely different
tails. 
\end{remark}
Further, we describe the decay of correlations along the interface.
Here again the horizontal scale $H^2$ plays an important role.

\begin{theorem}\label{thm:DecayOfCorrelations}\sl
Let $H=\Hone(\lambda)$ be as defined in \eqref{eq:Hgam.Def}. There
exist positive constants $C$, $c$, and $\lambda_0$ such
that for every $\lambda\in(0,\lambda_0]$ and all $i$, $j\in(0,N)$
we have
\[
\Cov(X_i,X_j)\le CH^{5/2}\exp\bl\{-c\,|i-j|\,H^{-2}\br\}.
\]
\end{theorem}

Thinking of $N$ as of time parameter, our system can be described as a
Markov chain on the positive half-line with certain attraction to the
origin. Being positive recurrent, its distribution $\mu_N^\lambda$ at
time $N$ approaches its stationary distribution $\pi_\lambda$
exponentially fast on the horizontal time scale $H^2$:

\begin{theorem}\label{thm:Relaxation}\sl
Let $H=\Hone(\lambda)$ be as defined in \eqref{eq:Hgam.Def}. There
exist positive constants $C$, $c$, and $\lambda_0$ such
that for every $\lambda\in(0,\lambda_0]$ we have
\[
\bl\|\mu_N^\lambda-\pi_\lambda\br\|_{\sf TV}
\le CH^2\exp\bl\{-cN\,H^{-2}\br\},
\]
where $\|\punkt\|_{\sf TV}$ denotes the total variational distance
between the probability measures.
\end{theorem}
\begin{remark}\sl
The reader might wonder whether the appearance of the exponents $1/3$,
$2/3$ and $3/2$ in the case $V(x)=|x|$ hints at a relationship between
the critical behaviour of the model considered here and the
much-studied Tracy-Widom distribution. We do not have a precise
answer to
this question; however at a heuristic level, the appearance of the
same critical exponents can be understood by noticing the similarities
between our model and the multi-layer PNG model introduced
in~\cite{mPhS02}, whose relation with the Tracy-Widom
distribution has been studied in the latter work. 
\end{remark}
\section{Proof of Theorem~\ref{thm:L1Bounds}}

\subsection{A basic comparison of partition functions}
\begin{lemma}
\label{lem_lowerboundPF}\sl
For any fixed $\rho>0$ and $\Hone(\lambda)$ defined as in
\eqref{eq:Hgam.Def}, put 
\begin{equation}\label{eq:Hdef}
H=\rho\Hone(\lambda).
\end{equation}
Then there exist positive constants $\lambda_0$, $c$, and $C$ such
that, for any $0<\lambda\le\lambda_0$, every $N\ge H^2$, and all
boundary conditions $0\le a,b\le H$, one has
\begin{equation}\label{eq:Basic.Comparison}
c\,e^{-CNH^{-2}}e^{-\lambda V(2H)N}
\ZNpo{a,b}\leq\ZNpl{a,b}\leq \ZNpo{a,b}\,.
\end{equation}
\end{lemma}

\begin{remark}\sl
Observe that
\begin{equation}
\label{eq:RemoveLambda}
\frac{\ZNpl{a,b}}{\ZNpo{a,b}}
=\ENpo{a,b}\Bl[\exp\Bl\{-\lambda\sum_{j=1}^{N-1}V\bl(X_j\br)\Br\}\Br]
\end{equation}
and thus the lemma states that for any $\lambda>0$ the exponential
moment of the functional 
\[
\BbbV(\BbbX)=\sum_{j=1}^{N-1}V\bl(X_j\br)
\] 
decays no faster
than exponentially in $N$ indicating that the typical value of
$V\bl(X_j\br)$ (equivalently, the height $X_j$ of the interface) is
``bounded on average''. It is instructive to compare this property to
the asymptotics 
\[
\ENpo{a,b}\Bl[N^{-3/2}\,\sum_{j=1}^NX_j\Br]\to \cst>0
\]
following from \cite{eB76,gL84}.
\end{remark}

As a straightforward application of the lemma above we get the
following simple but quite useful fact.

\begin{corollary}
\label{cor_freeRWcomparison}\sl
Under the conditions of Lemma~\ref{lem_lowerboundPF}, for any
collection $\calA\subset\INp$ of trajectories we have:
\begin{equation}
\label{eq:RW.UpperBound}
\PNpl{a,b}\bl(\calA\,\br)\le
\frac1c\exp\Bl\{CNH^{-2}+\lambda V(2H)N\Br\}\PNpo{a,b}\bl(\calA\,\br).
\end{equation}
Similarly, for any $\veps>0$, we get
\begin{equation}
\label{eq:RW.UpperBound.Area}
\begin{split}
\PNpl{a,b}\Bl(\calA\,;&\sum_{i=1}^N X_i\ge \veps HN\Br)
\\&\le\frac1c\exp\Bl\{CNH^{-2}+\lambda V(2H)N-\lambda V(\veps H)N\Br\}
\PNpo{a,b}\bl(\calA\,\br)
\end{split}
\end{equation}
and for any constants $A$, $B>0$,
\begin{equation}
\label{eq:RW.UpperBound.V}
\begin{split}
\PNpl{a,b}\Bl(\calA\,;&\BbbV(\BbbX)\ge AV(BH)N\Br)
\\&\le\frac1c\exp\Bl\{CNH^{-2}+\lambda V(2H)N-\lambda AV(BH)N\Br\}
\PNpo{a,b}\bl(\calA\,\br).
\end{split}
\end{equation}
\end{corollary}

\begin{proof}
Indeed, \eqref{eq:RW.UpperBound.Area} follows immediately from
convexity of $V(\punkt)$, as then
\[
\BbbV(\BbbX)=\sum_{j=1}^NV\bl(X_j\br)\ge NV\Bl(\sum_{j=1}^NX_j/N\Br)
\ge NV\bl(\veps H\br).
\]
Other inequalities are obvious.
\end{proof}

\begin{proof}[Proof of Lemma~\ref{lem_lowerboundPF}]
As the upper bound is obvious, one only needs to check the left
inequality. 

\begin{figure}[ht]
\bigskip
\centerline{\includegraphics[height=20mm]{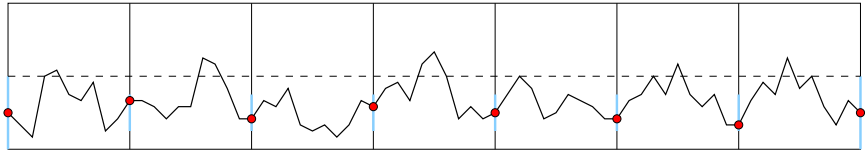}}
\bigskip
\caption{The event in the proof of Lemma~\ref{lem_lowerboundPF}.}
\label{fig_lowerboundPF}
\end{figure}

We use a renormalisation argument. With $H$ defined as in
\eqref{eq:Hdef},  take positive $\veps$ 
small enough to satisfy $\veps\in(0,1/4]$ and cut
each trajectory of our RW into blocks of length
\footnote{
both $H$ and $\Delta$ are assumed to be integer 
}
$\Delta=\veps H^2$;
this generates $\neps=\bl[N/\veps H^2\br]\ge4$ such blocks
(if there is a shorter piece left, we attach it to the last block). 
Further, denote $n_1=\Delta$, $n_2=(\neps-1)\Delta$ and 
consider the events 
\begin{gather*}
\calA=\Bl\{\every j=2,\dots,\neps-2:
\frac14H\le X_{j\Delta}\le\frac34H\Br\}\,,
\\[1ex]
\calB=\Bl\{\every j=n_1+1,\dots,n_2-1:0\le X_j\le 2H\Br\}\,,
\\[1ex]
\calC_1=\Bl\{\every j=1,\dots,n_1-1:0\le X_j\le 2H\Br\}
\cap\Bl\{\frac14H\le X_{n_1}\le\frac34H\Br\}\,,
\\[1ex]
\calC_2=\Bl\{\every j=n_2+1,\dots,N-1:0\le X_j\le 2H\Br\}
\cap\Bl\{\frac14H\le X_{n_2}\le\frac34H\Br\}\,.
\end{gather*}
Using \eqref{eq:RemoveLambda}, we immediately get
\begin{equation}\label{eq:ExpLowerBound}
\frac{\ZNpl{a,b}}{\ZNpo{a,b}}\ge\exp\Bl\{-\lambda V\bl(2H\br)N\Br\}\,
\PNpo{a,b}\bl(\calA\cap\calB\cap\calC_1\cap\calC_2\br)
\end{equation}
and it remains to bound below the probability of
$\calA\cap\calB\cap\calC_1\cap\calC_2$.  

First, observe that in view of the Donsker invariance principle we have 
\begin{gather*}
\min_{H/4\le x\le3H/4}\,\Pp{\Delta}x\bl(H/4\le X_\Delta\le3H/4\br)\ge C_1\,,
\\[1ex]
\min_{H/4\le c,d\le3H/4}
\PDpo{c,d}\bl(\every j=1,\dots,\Delta-1:0\le X_j\le 2H\br)\ge C_2\,.
\end{gather*}
By the Markov property we thus estimate
\[
\PNpo{a,b}\bl(\,\calA\cap\calB\mid\calC_1\cap\calC_2 \,\br)
\ge\bl(C_1C_2\br)^{\neps-2}\ge c\exp\bl\{-C\neps\br\}
\]
with some $c$, $C>0$.
On the other hand, by (conditional) independence of $\calC_1$ and
$\calC_2$ it is sufficient to bound below the probabilities
$\PNpo{a,b}\bl(\,\calC_1\,\br)$ and $\PNpo{a,b}\bl(\,\calC_2\,\br)$.
We shall estimate the former, the latter will follow in a similar way.

Combining the argument above with the conditional invariance principle
due to Bolthausen~\cite{eB76} we get 
\begin{equation}\label{eq:FromDonskerAndBolthausen}
\min_{a\in[H/4,H]\cup\{0\}}\Pp{\Delta}a\bl(\,\calC_1\,\br)\ge C_3>0
\end{equation}
with some constant $C_3=C_3(\veps)$, uniformly in $\lambda>0$ small
enough. Observe that the last bound holds also for any
$\widetilde\Delta$, $\Delta/2\le\widetilde\Delta\le\Delta$, and
perhaps smaller $C_3>0$. Consequently, it is enough to show that for
some constant $C_4=C_4(\veps)>0$ we have, uniformly in sufficiently
small $\lambda>0$, 
\begin{equation}\label{eq:StoppingTime}
\min_{0<a<H/4}\Pp{\Delta}a
\Bl(\,\min\Bl\{j\ge0:X_j\in\bl[H/4,H\br]\Br\}\le\frac\Delta4
\,\Br)\ge C_4>0
\end{equation}
as then immediately
\begin{equation}\label{eq:ComplementToBolthausen}
\min_{0< a<H/4}\Pp{\Delta}a\bl(\,\calC_1\,\br)\ge
C_3\,C_4>0\,.
\end{equation}

To check \eqref{eq:StoppingTime}, we fix an arbitrary integer $a$,
$0<a<H/4$, and consider two independent trajectories $X^a$ and $X^0$
distributed according to $\Pp{\Delta}a(\punkt)$ and
$\Pp{\Delta}0(\punkt)$ respectively. Let $\calD=\calD(k)$ denote the
``crossing event'' at $k$,
\[
\calD(k)=\Bl\{\,\every j=0,1,\dots,k-1: X^a_j>X^0_j
\quad\text{ \small and }\quad X^a_k\le X^0_k\,\Br\},
\]
and denote $\calD_\Delta=\cup_{k=1}^{n_3}\calD(k)$, $n_3=[\Delta/4]$. 
Our aim is to show that there exists a positive constant
$C_5=C_5(\veps)$ such that, uniformly in small enough $\lambda>0$,
one has (here and below, $\calD_\Delta^\comp$ stands for the
complement of $\calD_\Delta$)
\begin{equation}\label{eq:C1CouplingBound}
\Pp{\Delta}a\bl(\,\calC_1\bigm|\calD_\Delta\,)\ge C_5\,,
\qquad
\Pp{\Delta}a\bl(\,\calC_1\bigm|\calD_\Delta^\comp\,)\ge C_5\,.
\end{equation}
Then the target inequality
\[
\PNpo{a,b}\bl(\,\calC_1\,\br)\ge C_6(\veps)>0
\]
follows immediately. 

The key observation towards \eqref{eq:C1CouplingBound} is the
following. For an integer $x\ge0$, the jump distribution $p_x(\punkt)$
of our random walk from $x$ is given by
\[
p_x(k)=\frac{p(k)}{\rmP(\xi\ge-x)}\ind{\{k\ge-x\}}\,,
\]
where $\xi$ is a random variable with the unconstrained jump
distribution $p(\punkt)$. 
Clearly, the mean $e_x$ and the variance $\sig^2_x$ of $p_x(\punkt)$
satisfy
\begin{equation}\label{eq:PxMoments}
\begin{gathered}
e_x=\sum_kkp_x(k)\searrow e_\infty\equiv0
\quad\text{ as }x\nearrow\infty \\[1ex]
\sig^2_x=\sum_kk^2p_x(k)-\bl(e_x\br)^2
\le\frac{\sig^2}{\rmP(\xi\ge-x)}
\end{gathered}
\end{equation}
with $\sig^2$ denoting the variance of $p(\punkt)$. 

Now, suppose that the crossing event $\calD(k)$,
$k\in\{1,2,\dots,n_3\}$,  takes place. In view of
\eqref{eq:PxMoments}, with positive probability we have
$X^0_k-X^a_k\le M$, where the constant $M$ is independent of
$\lambda$. Thanks to the analogue of the aperiodicity property
\eqref{eq:aperiodic} for the distributions $p_x(\punkt)$ and
$p_y(\punkt)$, where $x=X^a_k$ and $y=X^0_k$ (with the same lower
bound \eqref{eq:aperiodic} for all $x,y\ge0$), two independent
trajectories started at $x$ and $y$ meet with positive probability
within $AM$ steps. Thus, the first inequality in
\eqref{eq:C1CouplingBound} follows immediately from the standard
independent coupling and the properties of $\Pp{\Delta}0(\punkt)$
mentioned above.

If $\calD_\Delta^\comp$ takes place, we have $X^a_j>X^0_j$ for all
$j=0,1,\dots,n_3$. Consequently, with positive probability the
stopping time 
\[
\tau=\min\bl\{j\ge0:X^a_j\ge H/4\br\}
\]
satisfies $\tau\le n_3$. Then the finite variance argument used
above implies that, with positive probability, we have 
$H/4\le X^a_\tau\le H$ and the second inequality in
\eqref{eq:C1CouplingBound} follows from a straightforward
generalization of~\eqref{eq:FromDonskerAndBolthausen}. 
\end{proof}

A literal repetition of the argument above implies also
the following result:

\begin{corollary}\label{cor:RemovingLambda}\sl
For positive $\rho$ and $\lambda$, put $H=2\rho\Hone(\lambda)$ and
define the event
\begin{equation}
\label{eq:EventBDef}
\calB=\calB_{H,N}=\Bl\{\every j=1,\dots,N:0\le X_j\le 2H\Br\}\,.
\end{equation}
Then, for any $\rho$, $\eta>0$, there exist positive constants
$\lambda_0$, $c$, and $C$ such that for all $0<\lambda\le\lambda_0$,
$N\ge\eta H^2$ and $0\le a,b\le H$, we have
\begin{equation}
\label{eq:EventABound}
\PNpo{a,b}\bl(\,\calB\,\br)\ge c\exp\bl\{-CNH^{-2}\br\}.
\end{equation}
Moreover, for any other event $\calA$ with 
$\PNpo{a,b}\bl(\,\calA\mid \calB\,\br)>0$, we get 
\begin{equation}
\label{eq:RemovingLambda}
c\,e^{-CNH^{-2}-2\lambda V(2H)N}
\le\frac{\PNpl{a,b}\bl(\,\calA\mid \calB\,\br)}
{\PNpo{a,b}\bl(\,\calA\mid \calB\,\br)}
\le c^{-1}\,e^{CNH^{-2}+2\lambda V(2H)N}\,,
\end{equation}
and thus, uniformly in bounded $NH^{-2}$, both conditional
probabilities are positive simultaneously.
\end{corollary}

\begin{proof}
Since $V(\punkt)$ is non-negative and monotone, \eqref{eq:EventBDef}
implies that, for any event $\calA$,
\[
e^{-\lambda V(2H)N}\ZNpo{a,b}\,\PNpo{a,b}\bl(\calA\calB\br)
\le \ZNpl{a,b}\PNpl{a,b}\,\bl(\calA\calB\br)
\le \ZNpo{a,b}\PNpo{a,b}\,\bl(\calA\calB\br)\,.
\]
The inequality \eqref{eq:RemovingLambda} now follows immediately from
Lemma~\ref{lem_lowerboundPF}. 
\end{proof}

\begin{remark}\label{rem:LambdaVsN}\sl
Although the importance of the scale $H=H_1$ (see~\eqref{eq:Hgam.Def})
should be clear from the proofs above, it is instructive to give
another motivation for the definition~\eqref{eq:Hgam.Def}.
Clearly, each interface under consideration can be naturally
decomposed into elementary excursions above the wall. Without external
field (ie., with $\lambda=0$) each such excursion of horizontal length
$l^2$ has typical height of order $l$. For $\lambda>0$ its
energetic price is of order at most $\lambda l^2V(l)$ and thus the
interaction with the field is negligible if $\lambda l^2V(l)\ll1$, that
is if $l=o(\Hone)$.
In other words, the presence of the field $\lambda$ is felt on the
(vertical) scale $\Hone$ or larger.
\end{remark}

\subsection{The upper bound}

The first half of Theorem~\ref{thm:L1Bounds}, namely
\[
\PNpl{}\Bl[N^{-1}\sum_{i=1}^N
     X_i \geq \delta^{-1}H\Br] 
\leq \frac1c\,e^{-C\delta^{-1}NH^{-2}},
\]
(with $H=\Hone(\lambda)$, see \eqref{eq:Hgam.Def})
follows directly from \eqref{eq:RW.UpperBound.Area} and the inequality
\[
V(\delta^{-1}H)\ge(2\delta)^{-1}V(2H)
\] 
valid for any $H\ge0$ and $0<2\delta\le1$:
\[
\PNpl{}\Bl[N^{-1}\sum_{i=1}^N X_i \geq \delta^{-1}H\Br]
\leq \frac1c\,e^{CNH^{-2}-((2\delta)^{-1}-1)\lambda V(2H)N}
\leq \frac1c\,e^{-C_\delta NH^{-2}};
\]
here $C_\delta=1/(4\delta)$ and $\delta$ is chosen small
enough to satisfy 
\[
0<\delta\le\min\Bl(1,\frac1{4(C+1)}\Br)\,.
\]
\qed

\subsection{The lower bound}\label{sec:LowerBound}

Our proof of the lower bound (with $H=\Hone(\lambda)$, see
\eqref{eq:Hgam.Def}), 
\begin{equation}
\label{eq:SecondHalf}
\PNpl{}\Bl[N^{-1}\sum_{i=1}^N
     X_i \leq \delta H\Br] 
\leq c\,e^{-C\delta^{-2}NH^{-2}},
\end{equation}
is based upon a certain renormalisation procedure. 
Namely, take $\lambda>0$ small enough and $\veps>0$ to be chosen later
(assuming, without loss of generality, that $\veps^2\in(0,1/12)$
and $\veps^2H^2$ is
an integer number larger than~$1$) and split every trajectory of the
random walk under consideration into pieces of length
$4\,\veps^2H^2$ to be called blocks; clearly, there are
exactly $\neps=\bl[N/(4\,\veps^2H^2)\br]$ such blocks (and, 
perhaps, an additional piece of shorter length). Next, we split each
block into four equal 
parts and use $\calI_m$, $m=1,\dots,4\neps$ to denote all obtained
sub-blocks. Further, we fix a small enough $\rho>0$ and say that the
trajectory under consideration is $\rho$-{\it high\/} in the $k$-th
block, if 
\[
\max_{j\in\calI_{4k+2}}X_j>\rho\,\veps H, \qquad
\max_{j\in\calI_{4k+4}}X_j>\rho\,\veps H \,.
\]
The main idea behind the argument below is as follows: for $\rho>0$
small enough, the number of $\rho$-high blocks in a typical trajectory
is of order $\neps$; however, a typical contribution of a $\rho$-high
block to the total area is of order at least
$\rho\,\veps^3H^3$; as a result, the typical area is bounded
below by a quantity of order
at least $\rho\,\veps^3H^3\neps\asymp\rho\,\veps HN$ and 
thus, for $\delta>0$ small enough, the event 
\[
\Ad\equiv\Bl\{\,\sum_{i=1}^NX_i \leq \delta HN\,\Br\}
\]
falls into the large deviation region for the distribution under
consideration. The target inequality \eqref{eq:SecondHalf} gives a
quantitative estimate for this to happen.

To start, we use \eqref{eq:RW.UpperBound} to remove the external field
$\lambda$, 
\[
\begin{split}
\PNpl{}\bl[\,\Ad\,\br] 
&\le c^{-1}\exp\Bl\{CNH^{-2}+\lambda V(2H)N\Br\}\PNpo{}\bl[\,\Ad\,\br]
\\[1ex]
&=c^{-1}\exp\Bl\{(C+1)NH^{-2}\Br\}\PNpo{}\bl[\,\Ad\,\br]\,.
\end{split}
\]
Now, conditionally on the configuration $\BbbX$ in the blocks
$\calI_{2m+1}$, $m=0,\dots,\neps-1$, the events
\[
\Bl\{\max_{j\in\calI_{2m}}X_j>\rho\,\veps H \Br\}
\]
are mutually independent. Moreover, 
a straightforward generalization of the argument used to estimate
below the probability of $\calC_1$ in Lemma~\ref{lem_lowerboundPF}
shows that 
\[
a_\rho\equiv\sup_m\sup_{0\le a,b\le\rho\,\veps H}\Ppo{\calI_{2m}}{a,b}
\Bl(\max_{j\in\calI_{2m}}X_j\le\rho\veps H\Br)\le1-\eta
\]
with some $\eta>0$, 
uniformly in $0<\rho\le\rho_0$ and all $\veps H=\veps\Hone(\lambda)$
large enough, implying that the events
\[
\Bl\{\text{ \small $k$-th block is $\rho$-high }\Br\}\,,
\]
which are also conditionally independent for fixed configuration in
$\calI_{4k+1}$, $k=0,\dots,\neps-1$, occur with probability at least
$(1-a_\rho)^2\ge\eta^2$. As a result, the number
$n_\rho$ of $\rho$-high blocks for a typical trajectory is not less
than $\eta^2\neps/2$. More precisely, since the events under
consideration are independent for individual blocks (conditionally on
every fixed configuration in between), the standard large deviation
bound implies 
\[
\PNpo{}\Bl(\,n_\rho<\frac{\eta^2}2\neps\Br)
\le \exp\bl\{\,-c\neps\,\br\}
=\exp\bl\{\,-c'\veps^{-2}NH^{-2}\,\br\}
\]
with some $c'=c'(\rho)>0$ not depending on $\veps$. Thus, taking
$\veps>0$ small enough, we obtain
\[
\PNpl{}\Bl(\,n_\rho<\frac\neps4\Br)
\le\exp\bl\{\,-c_1\veps^{-2}NH^{-2}\,\br\}.
\]

\begin{figure}[ht]
\centerline{\includegraphics[height=15mm]{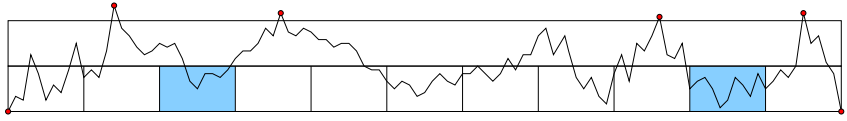}}
\caption{Two $\rho$-high blocks with oscillation}
\label{fig_hills}
\end{figure}

\indent From now on, we shall restrict ourselves to the trajectories
containing at least $\eta^2\neps/2$ blocks that are $\rho$-high.
We shall say that a $\rho$-high block oscillates if 
\[
\min_{j\in\calI_{4k+3}}X_j<\frac{\rho\,\veps}{2}H
\]
and observe that each $\rho$-high block without oscillation
contributes an amount at least $\rho\,\veps^3H^3/2$ to the total
area. Our final step of the proof consists in evaluating the typical
amount of oscillating blocks.

Define 
\[
l_k=\min\Bl\{j\in\calI_{4k+2}:X_j>\rho\,\veps H\Br\},\quad
r_k=\max\Bl\{j\in\calI_{4k+4}:X_j>\rho\,\veps H\Br\}
\]
and put $a_k=X_{l_k}$, $b_k=X_{r_k}$, $L=4k\veps^2H^2$, 
$R=4(k+1)\veps^2H^2$.
Then using essentially the same arguments as in the proof of
Lemma~\ref{lem_lowerboundPF}, we deduce that 
\[
\begin{split}
\min_{a,b}\,
\Ppo{(L,R)}{a,b}&\bl(\text{ \small $k$-th block does not oscillate }\br)
\\
&\ge\min\mathrm{P}_{r_k-l_k,+,0}^{\,a_k,b_k}
\Bl(\min_jX_j>\rho\,\veps H/2\Br)\ge 1-\bar a_\rho\,,
\end{split}
\]
where the bound $\bar a_\rho<1$ holds for any fixed $\rho$ small enough.
Arguing as before, we deduce that the number $\bar n_{\sf osc}$ of
blocks without oscillation satisfies the estimate
\[
\PNpo{}\Bl(\,\bar n_{\sf osc}<\frac{1-\bar a_\rho}2n_\rho\Br)
\le \exp\bl\{\,-c''n_\rho\,\br\}
=\exp\bl\{\,-c'''\veps^{-2}NH^{-2}\,\br\}
\]
with some $c'''=c'''(\rho)>0$ not depending on $\veps$. Again,
taking $\veps>0$ small enough, we obtain
\[
\PNpl{}\Bl(\,\bar n_{\sf osc}<\frac{1-\bar a_\rho}2n_\rho\Br)
\le \exp\bl\{\,-c_2\veps^{-2}NH^{-2}\,\br\}.
\]

However, on the complementary event we get
\[
\bar n_{\sf osc}\ge\frac{1-\bar a_\rho}2\,n_\rho
\ge \frac{1-\bar a_\rho}4\,\eta^2\,\neps
\]
for all $\rho>0$ small enough and thus the inequality
\[
\sum_{i=1}^NX_i \ge\frac{\rho\,\veps^3H^3}2\bar n_{\sf osc}
>\frac{\rho(1-\bar a_\rho)}{64}\,\eta^2\veps HN
\]
renders the event $\Ad$ impossible for 
$\delta=\rho(1-\bar a_\rho)\eta^2\veps/64$.
As a result, for such $\delta>0$ we get
\[
\begin{split}
\PNpl{}\bl[\,\Ad\,\br] 
&\le\PNpl{}\Bl(\,n_\rho<\frac{\eta^2\neps}2\Br)
+\PNpl{}\Bl(\,\bar n_{\sf osc}<\frac{1-\bar a_\rho}2n_\rho\Br)
\\[1ex]
&\le 2\exp\bl\{\,-c_3\delta^{-2}NH^{-2}\,\br\}
\end{split}
\]
with some $c_3=c_3(\rho)>0$.
\qed

\begin{remark}\sl
Obviously, the obtained lower $L^1$-bound on the total area 
implies immediately a simple lower bound for the height of the maximum
of interfaces: 
\begin{equation}
\label{eq:MaxLowerBound}
\begin{split}
\PNpl{}\Bl(\max_{1\le k\le N}X_k\le\delta\Hone \Br)
&\le\PNpl{}\Bl(N^{-1}\sum_{i=1}^N X_i \le \delta\Hone \Br)\\
&\le c\,\exp\Bl\{-C\delta^{-2}\lambda^{2/3}\, N\Br\}.
\end{split}
\end{equation}
We shall obtain a complementary bound after a more detailed analysis
of the interfaces.
\end{remark}

\section{Proof of Theorem~\ref{thm:PointwiseBounds}}

We treat the lower and the upper bounds in \eqref{eq:PointwiseBounds}
separately, the latter being based upon the following apriori estimates.

\subsection{Two refinements of the basic comparison lemma}

The following version of Lemma~\ref{lem_lowerboundPF} gives a better
bound than \eqref{eq:Basic.Comparison} for large values of $\rho$,
$\rho\ge\rho_0(\eta)>0$. 
With $H$ defined as in \eqref{eq:Hdef} and $\eta\in(0,1/2)$, we put
\footnote{
assuming $\Htil$ and $\Deltil$ to be integer.
}
\begin{equation}\label{eq:Htil}
\Htil=(1-2\eta)\rho\Hone(\lambda), 
\qquad \Deltil=\veps\Htil^2.
\end{equation}

\begin{lemma}\label{lem:LowerBoundRefined}\sl
Let $\rho$, $H$, $c$ and $C$ be as in Lemma~\ref{lem_lowerboundPF}. 
There exists $\lambda_0>0$ such that for any $\eta\in(0,1/2)$ and
$\zeta\in(0,1/2)$ there is
a constant $\tilde c>0$ such that for any $0<\lambda\le\lambda_0$,
every $N\ge H^2/(2\zeta)$, 
and all boundary conditions $0\le a,b\le H$, one has  
\begin{equation}\label{eq:Basic.Comparison.largeRho}
\tilde cc\,\exp\Bl\{-\frac{CN}{\Htil^2}-\lambda\Bl[\zeta V(2H)
+(1-\zeta)V(2\Htil)\Br]N\Br\}\le\frac{\ZNpl{a,b}}{\ZNpo{a,b}}
\le1\,.
\end{equation}
\end{lemma}

\begin{proof}
Let $H=\rho\Hone(\lambda)$, $\veps\in(0,1/4]$ and $\Delta=\veps H^2$
be as in the proof of Lemma~\ref{lem_lowerboundPF}. Similarly, for
$\eta\in(0,1/2)$ we define
\begin{equation}\label{eq:TildeParam}
\netil=\Bl[\frac{N-2\veps H^2}{\Deltil}\Br]
\equiv\Bl[\frac{N-2\veps H^2}{\veps(1-2\eta)^2\,H^2}\Br]\ge2\,.
\end{equation}
Further, let
\begin{gather*}
J_0=\Bl\{\Delta+j\Deltil: j=0,1,\dots,\netil-1\Br\}
\cup\bl\{N-\Delta\br\}
\\[1ex]
J_1=\Bl\{\Delta,\Delta+1,\dots,N-\Delta\Br\}, \qquad
J_2=\Bl\{1,\dots,N\Br\}\setminus J_1
\end{gather*}
\begin{figure}[ht]
\bigskip
\centerline{\includegraphics[height=25mm]{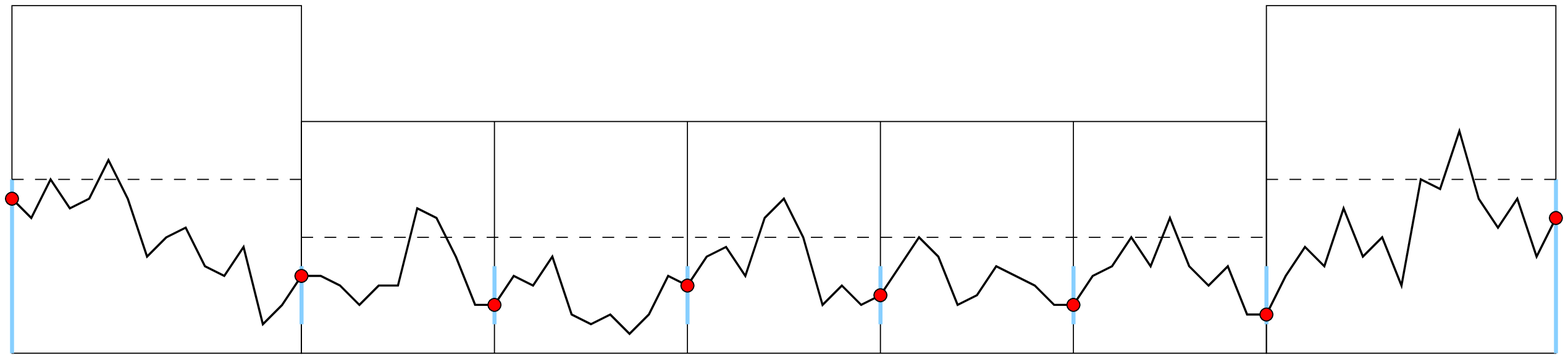}}
\bigskip
\caption{Renormalization scheme in Lemma~\ref{lem:LowerBoundRefined}.}
\label{fig:Renormalization}
\end{figure}
and introduce the events (see Fig.~\ref{fig:Renormalization}):
\begin{gather*}
\calA=\Bl\{\every j\in J_0:
\frac14\Htil\le X_{j}\le\frac34\Htil\Br\},
\\[1ex]
\calB=\Bl\{\every j\in J_1:0\le X_j\le 2\Htil\Br\}
\cap\Bl\{\every j\in J_2:0\le X_j\le 2H\Br\}\,.
\end{gather*}
For trajectories belonging to $\calA\cap \calB$ we have
\[
\BbbV(\BbbX)=\sum_{j=1}^{N-1}V\bl(X_j\br)
\le2\veps H^2V(2H)+(N-2\veps H^2)V(2\Htil)
\]
and therefore, denoting $\zeta=2\veps H^2/N$ we get 
(cf.~\eqref{eq:ExpLowerBound})
\[
\frac{\ZNpl{a,b}}{\ZNpo{a,b}}
\ge\exp\Bl\{-\lambda\Bl[\zeta V(2H)+(1-\zeta)V(2\Htil)\Br]N\Br\}\,
\PNpo{a,b}\bl(\calA\cap \calB\br).
\]
Moreover, using the scaling assumption \eqref{eq:TildeParam} and
arguing as in the proof of Lemma~\ref{lem_lowerboundPF} we get
\[
\PNpo{a,b}\bl(\calA\cap \calB\br)\ge ce^{-2\widetilde C\netil}
\ge c\,e^{-CN\Htil^{-2}}
\]
with  perhaps slightly smaller constant $\lambda_0>0$. 
\end{proof}

Next, we present a short-droplet analogue of the previous lemma.

\begin{lemma}\label{lem:ShortDropletLowerBound}\sl
Let $H=\rho\Hone(\lambda)$ and $N\le K\Hone^2(\lambda)$. 
There exist positive constants $\zeta$, $\lambda_0$ and $C$ such that
uniformly in $K/\rho^2<\zeta$, in $\lambda\in(0,\lambda_0]$ and all
boundary conditions $0\le a,b\le H$ one has
\[
C\ZNpo{a,b}\le\ZNpl{a,b}\le\ZNpo{a,b}\,.
\]
\end{lemma}

Our argument is based upon  the following small droplet bound to be
verified in Appendix~\ref{sec:SmallDropletBound} below.

\begin{lemma}\label{lem:SmallDropletBound}\sl
Let $S_0=0$, $S_k=\xi_1+\dots+\xi_k$, $k\ge1$, be the random walk
generated by a sequence $\xi_1$, $\xi_2$, \dots of i.i.d. random
variables such that $\bfE\xi=0$, $\bfE\xi^2=\sig^2<\infty$. 
Let $D>0$ be an arbitrary constant and, for any $m\ge1$, let $d_m$
satisfy $\bfP(S_m=d_m)>0$ and  $|d_m|\le D$.
Then there exists $\zeta>0$ such that
\begin{equation}\label{eq:SmallDropletBound}
\bfP\bl(\max_{0<k<m}S_k>M\mid S_m=d_m\br)\le\frac13,\qquad
\text{ as $M\to\infty$,}
\end{equation}
uniformly in $m/M^2\le\zeta$.
\end{lemma}

\begin{proof}[Proof of Lemma~\ref{lem:ShortDropletLowerBound}]
As in Lemma~\ref{lem_lowerboundPF}, our argument is based on
the bound (recall \eqref{eq:ExpLowerBound})
\[
\frac{\ZNpl{a,b}}{\ZNpo{a,b}}\ge
\exp\Bl\{-\lambda V\bl(5H\br)N\Br\}\,\PNpo{a,b}\bl(\calB\br)\,,
\]
where 
\[
\calB=\Bl\{\BbbX\in\INp:\max X_j\le5H\Br\};
\]
because of \eqref{eq:Hgam.Def}, \eqref{eq:Vregularity} and the
condition $N\le K\Hone^2(\lambda)$,
it remains to verify that the last probability is uniformly positive.
Let $\BbbX$ be an arbitrary trajectory from $\INp$ (recall
\eqref{eq:INp.Def}). Then, either it belongs to the set
\[
\calA_1=\Bl\{\BbbX:\every j=1,\dots,N-1:0\le X_j\le H\Br\}
\]
or there exists a non-empty set $[\frakl',\frakr']\subset(0,N)$ such that 
\begin{equation}\label{eq:LRprime.def}
\frakl'=\min\bl\{j>0: X_j>H\br\}\,,\quad
\frakr'=\max\bl\{j<N: X_j>H\br\}\,.
\end{equation}
We shall write $\Delta'=(\frakl',\frakr')$ and $|\Delta'|=\frakr'-\frakl'$.
Fix any $L\in(0,H)$ and denote
\begin{equation}\label{eq:A2.def}
\calA_2=\Bl\{\BbbX:X_{\frakl'}\le H+L,X_{\frakr'}\le H+L\Br\}\,.
\end{equation}
According to the bounded variance estimate~\eqref{eq:PxMoments} and
the Chebyshev inequality, one immediately gets
\begin{equation}\label{eq:L-Chebyshev.bound}
\PNpo{a,b}\bl(\BbbX\notin\calA_1\cup\calA_2\br)\le{2\sigtil^2}L^{-2}\,,
\end{equation}
where
\[
\sigtil^2=\max\Bl(\frac{\sig^2}{\rmP(\xi\ge0)},
\frac{\sig^2}{\rmP(\xi\le0)}\Br)\,.
\]
Taking $L$ sufficiently large to have 
$\PNpo{a,b}\bl(\BbbX\in\calA_1\cup\calA_2\br)\ge1/2$, we shall
restrict ourselves to trajectories $\BbbX$ belonging to
$\calA_1\cup\calA_2$ only (see Fig.~\ref{fig:shortdroplets}).
\begin{figure}[htb]
\bigskip
\centerline{\includegraphics[height=32mm]{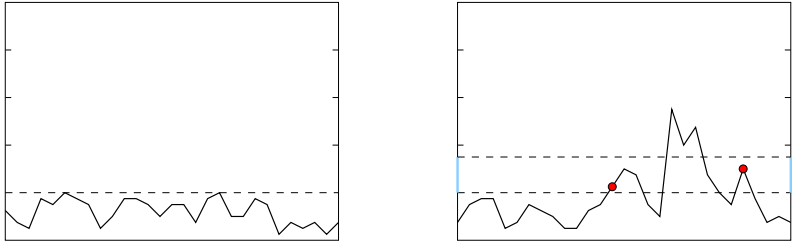}}
\bigskip
\caption{Trajectories from $\calA_1$ (left) and $\calA_2$ (right);
black dots correspond to the decomposition in \eqref{eq:notBonA2}.}
\label{fig:shortdroplets}
\end{figure}

As $\calA_1\subset\calB$, it remains to show that 
\begin{equation}\label{eq:BA2.bound}
\PNpo{a,b}\bl(\,\calB\mid\calA_2\,\br)\ge c
\end{equation}
for some constant $c>0$. Indeed, once \eqref{eq:BA2.bound} is
verified, we immediately get
\[
\begin{split}
\PNpo{a,b}\bl(\,\calB\,\br)
&\ge\PNpo{a,b}\bl(\,\calB\mid\calA_2\,\br)\PNpo{a,b}\bl(\calA_2\br)
+\PNpo{a,b}\bl(\calA_1\br)\\[1ex]
&\ge c\, \PNpo{a,b}\bl(\calA_1\cup\calA_2\br)\ge c/2
\end{split}
\]
and therefore
\[
\frac{\ZNpl{a,b}}{\ZNpo{a,b}}\ge
\exp\Bl\{-\lambda V\bl(5H\br)N\Br\}\,c/2\ge\tilde c
\]
uniformly in such $N$ and $\{a,b\}\subset[0,H]$.

To prove \eqref{eq:BA2.bound}, we rewrite
\begin{equation}\label{eq:notBonA2}
\begin{split}
\PNpo{a,b}\bl(\,\overline\calB\mid\calA_2\,\br)
&=\sum_{l',r'}\,\PNpo{a,b}
\bl(\frakl'=l',\frakr'=r'\mid \calA_2\br)
\\[1ex]
&\hphantom{=\sum_{l',r'}\,}\times\sum_{H\le a',b'\le H+L}
\Ppo{(l',r')}{\,a',b'}\bl(\,\overline\calB\,\br)
\\[1ex]
&\hphantom{=\sum_{l',r'}\,}
\times\PNpo{a,b}\bl(X_{l'}=a',X_{r'}=b'\mid\frakl'=l',\frakr'=r'\br)
\,.
\end{split}
\end{equation}
However, 
\begin{equation}\label{eq:CondBound}
\Ppo{(l',r')}{\,a',b'}\bl(\,\overline\calB\,\br)
\le\frac
{\mathrm{P}_{|\Delta'|}(\max S_j>3H\mid S_{|\Delta'|}=b'-a')}
{1-\mathrm{P}_{|\Delta'|}(\min S_j<-H\mid S_{|\Delta'|}=b'-a')}\,,
\end{equation}
where $\mathrm{P}_{|\Delta'|}$ refers to the distribution of
$|\Delta'|$-step unconstrained random walk with the step distribution
$p(\punkt)$, recall \eqref{eq:PDpl.Def}.
Finally, using the small droplet bound \eqref{eq:SmallDropletBound} and
taking $|\Delta'|/H^2\le N/H^2$ sufficiently small, we can make the
RHS above smaller than $1-c$. This finishes the proof. 
\end{proof}

\subsection{The upper bound}
\label{subsec:UpperBoundThm1.2}

We turn now to the proof of the upper bound in
Theorem~\ref{thm:PointwiseBounds}. Recall that due to the assumption
\eqref{eq:Vregularity} the function $V(\punkt)$ does not grow too fast
at infinity.

For $\rho>0$ and $\Hone=\Hone(\lambda)$, our canonical scale from
\eqref{eq:Hgam.Def}, define (cf.~\eqref{eq:DeltaLR.def})  
\begin{gather}
\label{eq:lrDef}
\begin{split}
&\frakl_\rho(M)=\max\Bl\{j<M: X_j\le\rho\Hone\Br\},\quad
\frakr_\rho(M)=\min\Bl\{j>M: X_j\le\rho\Hone\Br\}, \\[1ex]
&\hphantom{\frakl_\rho(M)}
\Delrho(M)=\bl\{\frakl_\rho+1,\frakl_\rho+2,\dots,\frakr_\rho-1\br\}, 
\qquad|\Delrho|=\frakr_\rho-\frakl_\rho-1
\end{split}
\end{gather}
and, for any integer interval $\Delta$,
\begin{equation}\label{eq:ArhoDef}
\calA_\rho(\Delta)=\Bl\{\every j\in\Delta, X_j>\rho\Hone\Br\}\,.
\end{equation}
Then, with $\Delta$ and $T$ being (large) natural numbers to be chosen
later, we get (for $N_2=[N/2]$ and $\Delrho=\Delrho(N_2)$)
\begin{equation}\label{eq:HeightUpperBound}
\begin{split}
\PNpl{}\Bl(X_{N_2}>T\Hone\Br)
&\le\PNpl{}\Bl(X_{N_2}>T\Hone\mid |\Delrho|<\Delta\Br)
\\[1ex]
&+\PNpl{}\bl(|\Delrho|\ge\Delta\br)\,.
\end{split}
\end{equation}

To estimate the length of the droplet, rewrite
\begin{equation}\label{eq:summation}
\begin{split}
\PNpl{}\bl(|\Delrho|\ge\Delta\br)
&=\sum_{l,r}\PNpl{}\bl(\Delrho=(l,r)\br)
\\[1ex]
&=\sum_{l,r}\sum_{0\le a,b\le\rho\Hone}\PNpl{}\bl(X_l=a,X_r=b\br)
\\[1ex]
&\hphantom{=\sum_{l,r}\sum_{0\le a,b\le H}}
\times\Ppl{(l,r) }{\,a,b}\bl(\calA_\rho(l,r)\br)\,,
\end{split}
\end{equation}
where the first summation goes over all $l$, $r$ satisfying
\begin{equation}\label{eq:lrSummationRegion}
0\le l<N_2< r\le N, \qquad r-l-1\ge\Delta\,.
\end{equation}
Next, by convexity of $V(\punkt)$ and the bounded growth assumption
\eqref{eq:Vregularity},
\begin{gather*}
V(H)\equiv V(\rho\Hone)\ge\frac\rho2\,V(2\Hone)
=\frac\rho{2\lambda\Hone^2},
\\[1ex]
\zeta V(2H)+(1-\zeta)V(2\Htil)
\le\Bl[\zeta f(2)+2(1-\zeta)(1-2\eta)\Br]V(H)\le\frac12V(H)
\end{gather*}
where $\Htil$ is as in \eqref{eq:Htil} and the constants $\zeta$, 
$\eta$ are chosen via 
\[
\zeta=\frac1{4f(2)}, \qquad \eta=\frac7{16}.
\]
Further,  applying Lemma~\ref{lem:LowerBoundRefined}, we obtain
(cf. Corollary~\ref{cor_freeRWcomparison}) 
\begin{equation}\label{eq:DropletEnergyBound}
\begin{split}
\Ppl{(l,r) }{\,a,b}\bl(\calA_\rho(l,r)\br)
&\le C_1\exp\Bl\{\Br[\frac C{\rho^2(1-2\eta)^2}-
\frac{\lambda V(H)\Hone^2}2\Br]\frac{r-l}{\Hone^2}\Br\}
\\[1ex]
&\hphantom{\le C'\exp}
\times\Ppo{(l,r) }{\,a,b}\bl(\calA_\rho(l,r)\br)
\\[1ex]
&\le C_1\exp\Bl\{-\Br[1-\frac{256C}{\rho^3}\Br]
\frac{\rho\Delta}{4\Hone^2}\Br\}
\,\Ppo{(l,r)}{\,a,b}\bl(\calA_\rho(l,r)\br).
\end{split}
\end{equation}
With $\Delta$, $\alpha\in(0,1)$ and $T>0$ satisfying
\begin{equation}\label{eq:AlphaTCondition}
\Delta=\sqrt T\Hone^2,\qquad 
\alpha T=\rho\ge \rho_0=8\cdot C^{1/3}\,,
\end{equation}
where $C$ denotes the same constant as in
Lemma~\ref{lem:LowerBoundRefined}, the last bound reads 
\[
\Ppl{(l,r)}{\,a,b}\bl(\calA_\rho(l,r)\br)
\le C_1\exp\Bl\{-\frac{\alpha}8\, T^{3/2}\Br\}
\,\Ppo{(l,r)}{\,a,b}\bl(\calA_\rho(l,r)\br)\,.
\]
Inserting it into \eqref{eq:summation} we immediately get
\begin{equation}\label{eq:DropletBound}
\PNpl{}\bl(|\Delrho|\ge\Delta\br)
\le C_1\exp\Bl\{-\frac{\alpha}8\, T^{3/2}\Br\}\,.
\end{equation}

It remains to estimate the first term in \eqref{eq:HeightUpperBound}.
Conditioning on the endpoints of the droplet of interest,
we decompose
\begin{equation}\label{eq:DropletDecomposition}
\begin{split}
\PNpl{}\bl(&X_{N_2}>T\Hone\mid{} |\Delrho|<\Delta\br)
\\[1ex]
&\le\sum_{l,r}\PNpl{}
\bl(\frakl_\rho=l,\frakr_\rho=r\mid |\Delrho|<\Delta\br)
\\
&\hphantom{\le\sum_{l,r}}
\times\sum_{0\le a,b\le H}
\PNpl{}\bl(X_l=a,X_r=b\mid\frakl_\rho=l,\frakr_\rho=r\br)
\\[1ex]
&\hphantom{\le\sum_{l,r}}
\times\Ppl{(l,r)}{\,a,b}\bl(X_{N_2}>T\Hone\mid\calA_\rho(l,r)\br)\,,
\end{split}
\end{equation}
where the first summation goes over all $l$, $r$ satisfying
(cf.~\eqref{eq:lrSummationRegion}) 
\[
0\le l<N_2< r\le N, \qquad r-l-1<\Delta\,.
\]
To finish the proof of the lemma it remains to establish the following
inequality 
\begin{equation}\label{eq:DropletTargetBound}
\Ppl{(l,r)}{\,a,b}\bl(X_{N_2}>T\Hone\mid\calA_\rho(l,r)\br)
\le\frac1{C_2}e^{-{C_2}T^{3/2}}\,.
\end{equation}
Notice that taking $T_0$ large enough, we can achieve the bound
\[
\frac\Delta{(\rho\Hone)^2}
\le\frac{\sqrt T\Hone^2}{(\alpha T\Hone)^2}
=\frac1{\alpha^2T_0^{3/2}}\le\zeta
\]
for all $T\ge T_0$, where $\zeta$ is the same constant as in
Lemma~\ref{lem:ShortDropletLowerBound}, and thus can remove the field
$\lambda$ from our further considerations.

\bigskip
\begin{figure}[ht]
\bigskip
\centerline{\includegraphics[height=35mm]{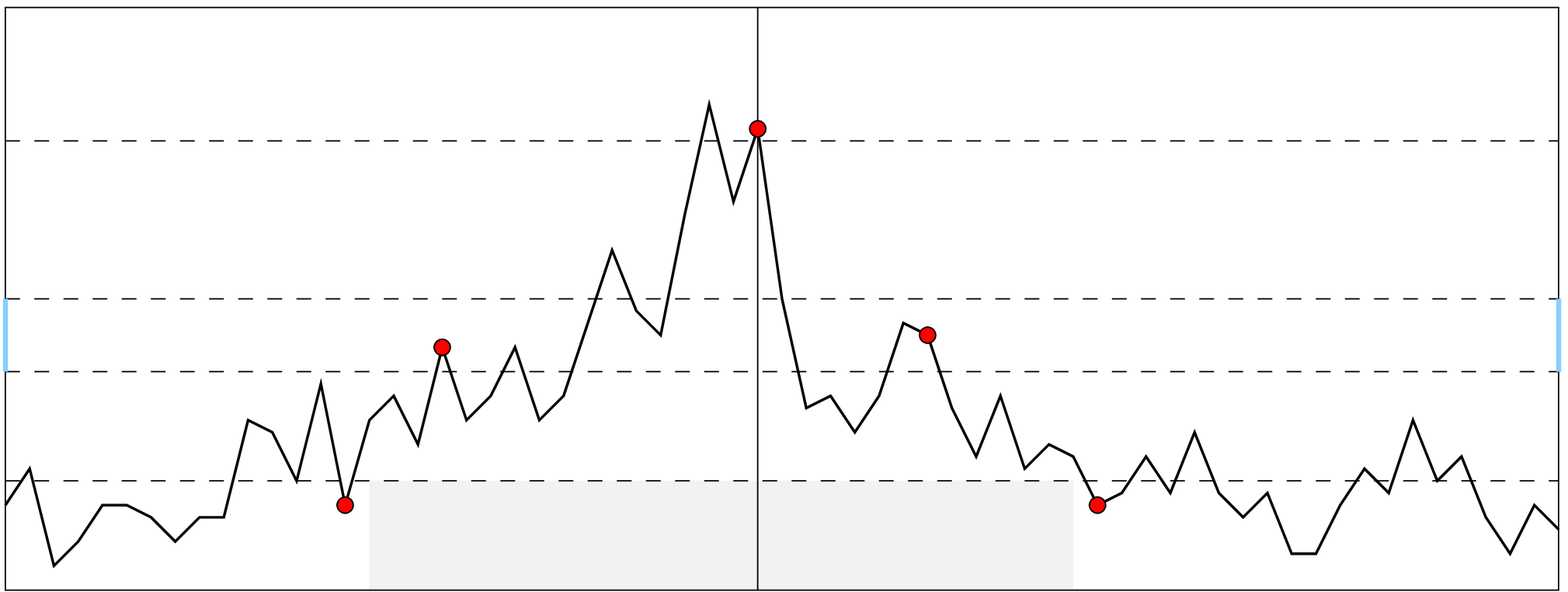}}
\bigskip
\caption{Decompositions \eqref{eq:DropletDecomposition} and (\ref{eq:notBonA2}$'$)}
\label{fig:highdroplet}
\end{figure}

To prove \eqref{eq:DropletTargetBound}, we shall proceed as in the
proof of Lemma~\ref{lem:ShortDropletLowerBound}. Namely, defining
$\Delta'\subset(l,r)$ and $\calA_2$ similarly to
\eqref{eq:LRprime.def} and  \eqref{eq:A2.def}, (see
Fig.~\ref{fig:highdroplet}), we bound above
\[
\Ppo{(l,r)}{\,a,b}\bl(X_{N_2}>T\Hone\mid\calA_\rho(l,r)\br)
\le\Ppo{(l,r)}{\,a,b}\bl(\overline{\calA}_2\br)+
\Ppo{(l,r)}{\,a,b}\bl(X_{N_2}>T\Hone\mid\calA_2\br)
\]
and decompose (cf.~\eqref{eq:notBonA2})
\begin{equation}\tag{\ref{eq:notBonA2}$'$}
\begin{split}
\Ppo{(l,r)}{\,a,b}&\bl(\,X_{N_2}>T\Hone\mid\calA_2\,\br)
=\sum_{l\le l'\le r'\le r}\,\Ppo{(l,r)}{\,a,b}
\bl(\frakl'=l',\frakr'=r'\mid \calA_2\br)
\\[1ex]
&\hphantom{\bl(\,X_{N_2}>T\Hone\mid}
\times\sum_{H\le a',b'\le H+L}
\Ppo{(l',r')}{\,a',b'}\bl(\,X_{N_2}>T\Hone\,\br)
\\[1ex]
&\hphantom{\bl(\,X_{N_2}>T\Hone\mid}
\times\Ppo{(l',r')}{\,a',b'}
\bl(X_{l'}=a',X_{r'}=b'\mid\frakl'=l',\frakr'=r'\br)
\,.
\end{split}
\end{equation}

We shall estimate the probability
$\Ppo{(l',r')}{\,a',b'}\bl(\,X_{N_2}>T\Hone\,\br)$ depending on the length
of the interval $\Delta'$. 
First, we observe that for any $L>0$ the
Donsker invariance principle gives (recall \eqref{eq:CondBound})
\[
\begin{split}
\Ppo{(l',r')}{\,a',b'}\bl(\,X_{N_2}>T\Hone\,\br)
&\le\frac{\mathrm{P}_{|\Delta'|}(\max S_j>(T-2\rho)\Hone
\mid S_{|\Delta'|}=b'-a')}
{1-\mathrm{P}_{|\Delta'|}(\min S_j<-\rho\Hone\mid S_{|\Delta'|}=b'-a')}
\\[1ex]
&\le C_3\exp\Bl\{-C_4\frac{(T-2\rho)^2}{\sqrt T}\Br\}
=C_3e^{-C_4(1-2\alpha)^2T^{3/2}}
\end{split}
\]
uniformly in $\Hone^{7/6}\le|\Delta'|\le\Delta=\sqrt T\Hone^2$, in
$\{a',b'\}\subset[H,H+L]$ and all $\lambda$ small enough.

On the other hand, for $|\Delta'|\le\Hone^{7/6}$, the conditional
Chebyshev inequality for the maximum (see
Lemma~\ref{lem:ConditionalChebyshevForMaximum} in the appendix below)
gives 
\[
\Ppo{(l',r')}{\,a',b'}\bl(\,X_{N_2}>T\Hone\,\br)
\le\frac{\cstL|\Delta'|^{3/2}(T-2\rho)^{-2}\Hone^{-2}}
{1-\cstL|\Delta'|\rho^{-2}\Hone^{-2}}
\le C_5\frac{\cstL}{T^2\Hone^{1/4}}\,.
\]

Now, combining the last two estimates with the bound
\eqref{eq:L-Chebyshev.bound}, we obtain:
\begin{equation}\label{eq:SmallDropletBoundTotal}
\Ppo{(l,r)}{\,a,b}\bl(\,X_{N_2}>T\Hone\mid\calA_2\,\br)
\le\frac{2\sigtil^2}{L^2}+C_3e^{-C_4(1-2\alpha)^2T^{3/2}}
+C_5\frac{\cstL}{T^2\Hone^{1/4}}\,.
\end{equation}
To finish the proof, we first take
$L=\exp\Bl\{\frac{C_4}2(1-2\alpha)^2T^{3/2}\Br\}$
and then $\lambda$ sufficiently small to make the last term smaller
than the second. With this choice \eqref{eq:SmallDropletBoundTotal}
reduces to \eqref{eq:DropletTargetBound}. 
\qed

An obvious generalization of the argument above gives also the
following bound.

\begin{corollary}
\label{cor_HighDropletBound}
\sl
Let $\rho$ and $\veps$ denote some (small) positive constants and let
the integers $l_0$, $r_0$, $M\in[0,N]$ be such that 
\[
l_0<M<r_0\qquad\text{ and }\qquad
\min\bl(|M-l_0|,|M-r_0|\br)\ge\veps\Hone^2
\] 
with $\Hone=\Hone(\lambda)$ being our canonical scale from
\eqref{eq:Hgam.Def}. Then there exists $c_1>0$ such that for any $T>0$
large enough and all $a_0$, $b_0\in[0,\rho\Hone]$ the inequality
\[
\Ppl{(l_0,r_0)}{a_0,b_0}
\bl(X_M>T\Hone\bigm|\calA_{\rho}(l_0,r_0)\br)\le e^{-c_1T^{3/2}}
\]
holds for all $\lambda\in(0,\lambda_0]$, where
$\lambda_0=\lambda_0(T)>0$. 
\end{corollary}

As a straightforward modification of the proofs above one can show
existence of moments of $X_{N_2}$ to be used below.

\begin{corollary}
\label{cor_HighMomentUpperBound}
\sl
There exist positive constants $K$ and $\lambda_0$ such that for all
$p$, $1<p<21/8$, we have
\[
\ENpl{}\bl(\,X_{N_2}\,\br)^{2p}\le C(p)\,\Hone^{2p+1}\,,
\]
uniformly in $\lambda\in(0,\lambda_0]$ and $N\ge K\Hone^2$.
\end{corollary}

\begin{proof}
Using the decomposition \eqref{eq:summation} with $\rho=T/4$ and
$\Delta=\Hone^2T^{-9/10}$ in \eqref{eq:AlphaTCondition}, we get the
following analogue of \eqref{eq:DropletBound}:
\[
\PNpl{}\bl(|\Delrho|\ge\Delta\br)
\le C_1\exp\Bl\{-\frac1{32}\,T^{1/10}\Br\}\,.
\]
Next, we use the decomposition (\ref{eq:notBonA2}$'$) and bound
above the height of the inner droplet 
\[
\Ppo{(l',r')}{\,a',b'}\bl(\,X_{N_2}>T\Hone\,\br)\,.
\]
As in the proof of Lemma~\ref{lem_ConditionalChebyshev} two cases
to be considered separately, $|r'-l'|\le m_0$ and
$m_0\le|r'-l'|\le\Delta$. Clearly, w.l.o.g.\ we may and shall assume
that $m_0$ is chosen large enough to satisfy
(cf.~\eqref{eq:LLTLowerBound}) 
\[
\bfP\bl(S_m=d_m\br)\ge\frac1{2e\sqrt{2\pi\sig^2m}}
\]
for all $m\ge m_0$.

Let $|r'-l'|\le m_0$. Combining \eqref{eq:FiniteDropletProbability}
and \eqref{eq:OnePointChebyshev}, we get
\[
\Ppo{(l',r')}{\,a',b'}\bl(\,X_{N_2}>T\Hone\,\br)
\le\frac{4\sig^4m_0{}^4}{p(m_0,D)}\bl(T\Hone)^{-4}\,.
\]
On the other hand, for $(l',r')$ satisfying $m_0\le|r'-l'|\le\delta$,
we apply Lemma~\ref{lem_OnePointConditionalChebyshev} to obtain
\[
\Ppo{(l',r')}{\,a',b'}\bl(\,X_{N_2}>T\Hone\,\br)
\le\frac{C\sig|r'-l'|^{5/2}}{(T\Hone)^4}
\]
with a numeric constant $C\le18$. As a result, for any $T>0$ we get
\[
\PNpl{}\bl(\,X_{N_2}>T\Hone\,\br)\le
C_1\exp\Bl\{-\frac1{32}\,T^{1/10}\Br\}
+\frac{C_2(m_0,D)}{(T\Hone)^4}+\frac{C_3\Hone}{T^{25/4}}\,,
\]
and therefore, for $1<p<21/8$,
\[
\begin{split}
\ENpl{}\bl(\,X_{N_2}\,\br)^{2p}&\le \Hone^{2p}
\Bl(1+\sum_{T>0} (T+1)^{2p}
\PNpl{}\bl[T\Hone \le X_i<(T+1)\Hone]\Br)\\[1ex]
&\le C_4(p)\,\Hone^{2p+1}\,.
\end{split}
\]
\end{proof}

\subsection{The lower bound}

Fix a (big) positive $T$ and an integer $\Delta=\sqrt{T}\Hone^2$ and
denote (recall the notation $N_2=[N/2]$)
\[
l=N_2-\Delta,\qquad r=N_2+\Delta\,.
\]
It follows from the argument in Sect.~\ref{subsec:UpperBoundThm1.2}
that for some constant $T_0>0$ not depending of $\lambda>0$ we have
\begin{equation}\label{eq:LowerBoundScale}
\PNpl{}\bl(\bl\{X_l,X_r\br\}\subset[0,T_0\Hone]\br)\ge\frac12
\end{equation}
(recall the running
assumption that we omit the boundary conditions $a=b=0$ from the
notation). We thus rewrite
\begin{equation}\label{eq:LowerBoundLocal}
\begin{split}
\PNpl{}\bl(X_{N_2}>T\Hone\br)\ge\sum_{0\le a,b\le T_0H}
&\PNpl{}\bl(X_l=a,X_r=b\br)
\\[1ex]
&\kern-13pt\times\Ppl{(l,r)}{\,a,b}\bl(X_{N_2}>T\Hone\br)\,.
\end{split}
\end{equation}
Now, define 
\[
l'=l+H^2,\qquad r'=r-H^2,\qquad X_{l'}=a',\qquad X_{r'}=b' 
\]
and estimate
\begin{equation}\label{eq:LowerBoundLocal2}
\begin{split}
\Ppl{(l,r)}{\,a,b}\bl(X_{N_2}>T\Hone\br)\ge\sum_{H/4\le a',b'\le3H/4}
&\Ppl{(l,r)}{\,a,b}\bl(X_{l'}=a',X_{r'}=b'\br)
\\[1ex]
&\kern-13pt\times\Ppl{(l',r')}{\,a',b'}\bl(X_{N_2}>T\Hone\br)\,.
\end{split}
\end{equation}
By the Donsker invariance principle and the estimates for the maximum
of the Brownian bridge, the last factor is bounded below by
\[
\frac1c\exp\bl\{-c(T\Hone)^2/(r'-l')\br\}
\ge\frac1{c'}\exp\bl\{-c'T^{3/2}\br\}
\]
uniformly in $a'$ and $b'$ under consideration, provided only
$\lambda>0$ is small enough. 
Next, a literal repetition of the proof of
\eqref{eq:FromDonskerAndBolthausen} and
\eqref{eq:ComplementToBolthausen} combined with the estimate
\eqref{eq:LowerBoundLocal2} gives
\begin{equation}\label{eq:LowerBoundLocal3}
\Ppl{(l,r)}{\,a,b}\bl(X_{N_2}>T\Hone\br)\ge
\frac1{c''}\exp\bl\{-c''T^{3/2}\br\}
\end{equation}
uniformly in $a$, $b$ from $[0,T_0\Hone]$ provided only $\lambda>0$ is
sufficiently small. 
The lower bound in \eqref{eq:PointwiseBounds} now follows from
\eqref{eq:LowerBoundScale}, \eqref{eq:LowerBoundLocal},
and~\eqref{eq:LowerBoundLocal3}. 
\qed

\section{Refined asymptotics}

\subsection{Quasirenewal structure}\label{sec:quasirenewal}

The importance of the scale $\Hone(\lambda)$ demonstrated in the
proofs of the previous sections is even more pronounced in the study
of the refined behaviour of the interfaces under consideration. The
aim of this section is to describe certain intrinsic renewal-type
structure of the random walks distributed via
\eqref{eq:INp.Def}--\eqref{eq:aperiodic} that manifests itself in the
diffusing scaling (i.e., $\Hone(\lambda)^2$ in the horizontal
direction and $\Hone(\lambda)$ in the vertical one).

For any $\rho>0$ and $\lambda>0$, let $S_\rho$ denote the horizontal
strip of width $4\rho\Hone$,
\[
S_\rho=\Bl\{(x,y)\in\BbbZ^2:y\in[0,4\rho\Hone]\Br\}.
\] 
In this section we shall 
establish certain quasirenewal property stating roughly that for all
$\lambda>0$ small enough the ``density'' of visits of the RW under
consideration to the strip $S_\rho$ is ``positive on the scale
$\Hone(\lambda)^2$''.

More precisely, for positive real $\veps$, $\rho$, $\lambda$ and
integer $K>0$ we split our trajectories into $K$-blocks $\calI_m$ of
length $K\veps\Hone(\lambda)^2$ (assuming w.l.o.g.\ $\veps\Hone(\lambda)^2$
to be integer) and introduce the random variables 
\begin{equation}
\label{eq:YLabels}
Y_m=Y_m^\rho
=\ind{\{\min_{j\in\calI_m}X_j>2\rho\Hone \}}\,.
\end{equation}
Let $n_K=\bl[N/(K\veps\Hone(\lambda)^2)\br]$ be the total number of
such blocks and let $n_Y$ be the total number of $K$-blocks labelled
by ones: 
\[
n_Y=\left|\bl\{m:Y_m=1\br\}\right|\,.
\]
Our first observation is that with high probability the total length
$n_YK\veps\Hone(\lambda)^2$ of such $K$-blocks can not be large:

\begin{lemma}\label{lem:nYBound}\sl
Let $f(\cdot)$ be defined as in \eqref{eq:Vregularity}, $\veps$ be a
positive constant, and $\alpha$ satisfy $\alpha\in(0,1)$.
For any $\rho>0$ there exist positive constants $\lambda_0$, and
$K_0$, depending on $\alpha$, $\veps$, and $\rho$ only,
such that for any $\lambda\in(0,\lambda_0)$, any
$K\ge K_0$ and any $N\ge3K_0\veps\Hone(\lambda)^2$ we have
\begin{equation}
\label{eq:nYBound}
\PNpl{a,b}\bl(n_Y\ge\alpha n_K\br)
\le\exp\Bl\{-\frac\alpha{8f(2/\rho)}N\Hone(\lambda)^{-2}\Br\}
\end{equation}
uniformly in $a$, $b\in(0,\rho\Hone )$.
\end{lemma}

\begin{proof}[Proof of Lemma~\ref{lem:nYBound}]
We use the blocking procedure described at the beginning of this
section where, given $\rho>0$ and $\veps>0$, the constant $K_0\ge8$ is
chosen large enough, see \eqref{eq:K0large} below.

The $Y$-labels defined in \eqref{eq:YLabels} with any $K\ge K_0$
introduce a $0-1$ encoding of each trajectory; using this encoding, we
split the $K$-blocks labelled by ones into maximal ``connected
components'' to be called $K$-clusters. Two neighbouring $K$-clusters
are called connected if they are separated by exactly one $K$-block
labelled by a zero. A $K$-cluster that is not connected to its
neighbours is called isolated. Our next goal is to show that for any
collection of $K$-clusters consisting of $n_Y$ $K$-blocks there is a
sub-collection of isolated $K$-clusters of total length at least
$\bl[n_Y/4\br]$. As soon as this is done, a simple reduction argument
will imply the target estimate \eqref{eq:nYBound}.

\begin{figure}[ht]
\bigskip
\centerline{\includegraphics[height=45mm]{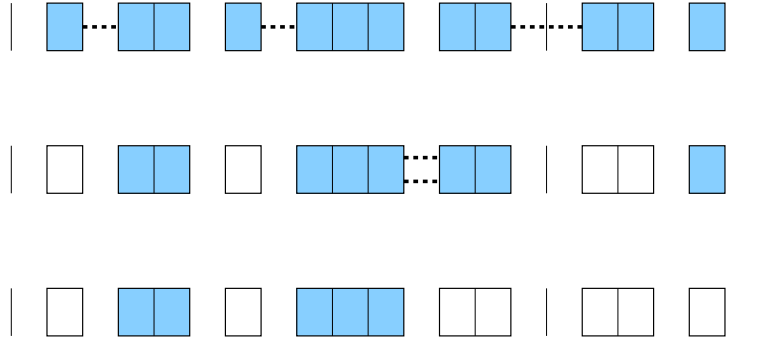}}
\bigskip
\caption{Two-step selection procedure: from $7$ clusters of total
length $12$ choose $2$ clusters of total length $5$; selected
$K$-clusters are shadowed} 
\label{fig:decimation}
\end{figure}

Our selection procedure consists of two steps. First, we split all
$K$-clusters into subsequent pairs of neighbouring clusters and from
each such pair we choose the longest cluster (or the left one if they
are of equal length). Observe that each chosen cluster is either
isolated or belongs to a pair of connected $K$-clusters. Next, we
split all isolated $K$-clusters into subsequent pairs
\footnote{
temporarily neglecting the very last $K$-cluster, if their total
number is odd
}
and from each pair (connected or isolated) we choose the longest
cluster (or the left one if they are of equal length). The obtained
collection (together with the temporarily neglected $K$-cluster, if
there was one) consists of isolated clusters which altogether contain
at least $n_Y/4$ of $K$-blocks, see Fig.~\ref{fig:decimation}.

Our second step relies upon a finer renormalisation, this time on the
integer scale $\veps\Hone(\lambda)^2$. We split our trajectory into
$\neps=\bl[N/\veps\Hone(\lambda)^2\bl]$ blocks $\calJ_l$ of length
$\veps\Hone(\lambda)^2$ each and similarly to \eqref{eq:YLabels}
introduce the labels
\[
Z_l=Z_l^\rho=\ind{\{\min_{j\in\calJ_l}X_j>2\rho\Hone \}}\,.
\]
Of course, the natural (inclusion) correspondence between
$\veps$-blocks $\calJ_l$ and $K$-blocks $\calI_m$,
\[
\calJ_l\subset\calI_m\,,
\]
has the following property: if $Y_m=1$
and $\calJ_l$ corresponds to $\calI_m$, then $Z_l=1$. As before, we
split all $\veps$-blocks labelled by ones into maximal connected
components to be called $\veps$-clusters. Clearly, as subsets of
$\bl\{1,\dots,N\br\}$, every $K$-cluster is included in the
corresponding $\veps$-cluster. Let $\calE$ be the collection of
$\veps$-clusters that correspond to the isolated $K$-clusters selected
by applying a procedure as in Fig.~\ref{fig:decimation}. 
The following two properties of the collection
$\calE$ will be important for our future application: 
1) every $\veps$-cluster from $\calE$ is bounded by two boundary
$\veps$-blocks labelled by zeroes; moreover, for different
$\veps$-clusters the boundary blocks are different;
2) the total length of $\veps$-clusters from $\calE$ is at least
$Kn_Y\veps\Hone(\lambda)^2/4$. 

Our next aim is to establish certain one-droplet estimate from which
the target bound \eqref{eq:nYBound} will follow immediately. Consider
any $\veps$-cluster from $\calE$ and denote its extremal $\veps$-blocks
(the first and the last one) by $\calJ_{m_1}$ and $\calJ_{m_2}$
respectively. Clearly, the length of this $\veps$-cluster is
$\neps^0\veps\Hone^2\equiv(m_2-m_1+1)\veps\Hone^2$, $\neps^0\ge K$. 
Further, define
\begin{gather*}
l=\max\Bl\{j\in\calI_{m_1-1}:X_j<2\rho\Hone \Br\},
\quad\quad X_l=a,\\
r=\min\Bl\{j\in\calI_{m_2+1}:X_j<2\rho\Hone \Br\},
\quad\quad X_r=b.
\end{gather*}
Similarly, let
\begin{equation}\label{eq:lr0Def}
l_0=m_1\veps\Hone^2,\quad X_{l_0}=a_0, \quad 
r_0=(m_2-1)\veps\Hone^2, \quad X_{r_0}=b_0\,.
\end{equation}
Using the notation $\calA_2=\calA_{2\rho}(l_0,r_0)$ 
(recall \eqref{eq:ArhoDef}), one
can bound above the partition function corresponding to
this droplet by
\begin{equation}\label{eq:DropPFDecompose}
\Zpl{(l,r)}{a,b}\bl(\calA_2\br)
=\sum_{a_0,b_0\ge2\rho\Hone}\Zpl{(l,l_0)}{a,a_0}
\Zpl{(l_0,r_0)}{a_0,b_0}\bl(\calA_2\br)
\Zpl{(r_0,r)}{b_0,b}\,.
\end{equation}
Clearly, the target inequality \eqref{eq:nYBound} follows immediately
from the one-droplet bound
\begin{equation}\label{eq:OneDropBound}
\sup_{\rho\Hone\le a_0,b_0\le2\rho\Hone}
\frac{\Zpl{(l_0,r_0)}{a_0,b_0}\bl(\calA_2\br)}{\Zpl{(l_0,r_0)}{a_0,b_0}}
\le \exp\Bl\{-\frac{\veps\neps^0}{2f(2/\rho)}\Br\}
\end{equation}
the lower bound on the total length of $\veps$-clusters from
$\calE$, provided only 
\[
K\veps\alpha>4f(2/\rho)\log2\,,
\]
to suppress the total number of $0-1$ encodings (that is bounded above
by $2^{n_K}$).

Our proof of \eqref{eq:OneDropBound} will be based upon the
decomposition \eqref{eq:DropPFDecompose} and the following two facts:
{\sl
\begin{quote}
F1) there is $\lambda_0=\lambda_0(\veps,\rho,...)>0$ for which: for
any $\eta>0$ there exists $T>0$ such that uniformly in
$\lambda\in(0,\lambda_0]$ one has (recall \eqref{eq:lr0Def})
\end{quote}
\begin{equation}\label{eq:RefinedFact1}
\Zpl{(l,r)}{a,b}\bl(\calA_2,a_0\le2T\Hone,b_0\le2T\Hone\br)
\ge(1-\eta)\Zpl{(l,r)}{a,b}\bl(\calA_2\br)
\end{equation}

\begin{quote}
F2) for $\veps$, $\rho$, $T$ and $\lambda_0$ as above there is a
finite constant $M>0$ such that, 
uniformly in $\lambda\in(0,\lambda_0]$,
\end{quote}
\begin{equation}\label{eq:RefinedFact2}
\sup
\frac{\Zpl{\Delta}{a,a_0}}{\Zpl{\Delta}{a,a_0'}}\le M
\end{equation}
\begin{quote}
with supremum taken over $\veps\Hone^2\le\Delta\le2\veps\Hone^2$, 
$\rho\Hone\le a\le2\rho\Hone$, and $2\rho\Hone\le
a_0\le(2T+\rho)\Hone$, $a_0'\equiv a_0-\rho\Hone$;
a similar estimate (with the same constant $M$) holds for the ratio
$\Zpl{\Delta}{a_0,a}/\Zpl{\Delta}{a_0',a}$.
\end{quote}
}

The inequality \eqref{eq:OneDropBound} follows easily from
\eqref{eq:RefinedFact1} and \eqref{eq:RefinedFact2}. Indeed, combining
\eqref{eq:DropPFDecompose} and \eqref{eq:RefinedFact1} we get
\[
\begin{split}
\Zpl{(l,r)}{a,b}\bl(\calA_2\br)&\le\frac1{1-\eta}
\Zpl{(l,r)}{a,b}\bl(\calA_2,a_0\le2T\Hone,b_0\le2T\Hone\br)
\\[1ex]
&\le\frac1{1-\eta}
\sum\Zpl{(l,l_0)}{a,a_0}\Zpl{(l_0,r_0)}{a_0,b_0}\bl(\calA_2\br)
\Zpl{(r_0,r)}{b_0,b}
\end{split}
\]
with the sum running over $2\rho\Hone\le a_0,b_0\le(2T+\rho)\Hone$.
Further, denoting $\calA_1=\calA_\rho(l_0,r_0)$, taking $\eta=1/2$, 
and using the estimate \eqref{eq:RefinedFact2}, the convexity of the
function $V(\,\cdot\,)$ and the reduction of the central part of the
droplet as in Fig.~\ref{fig_reduction}, we bound the last expression
by 
\begin{figure}[ht]
\bigskip
\centerline{\includegraphics[height=75mm]{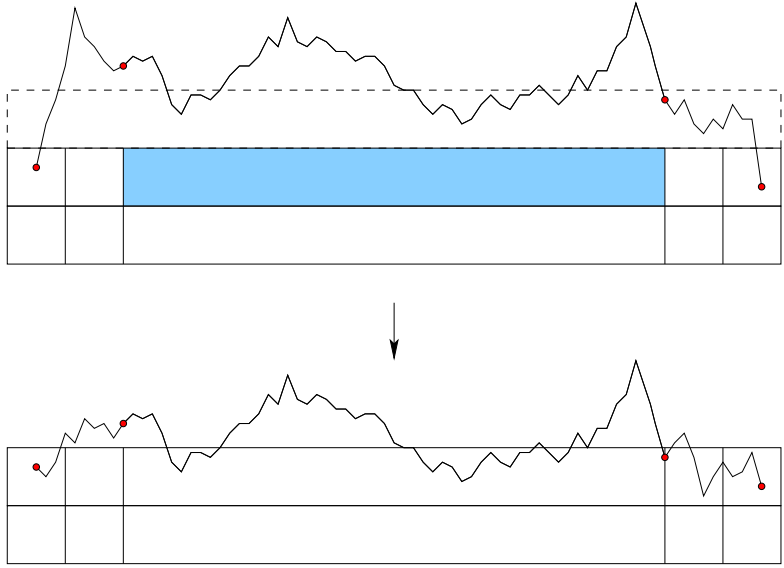}}
\bigskip
\caption{Reduction of an $\veps$-cluster}
\label{fig_reduction}
\end{figure}
\begin{gather*}
2M^2e^{-\lambda\neps^0\veps\Hone^2V(\rho\Hone)}
\sum\Zpl{(l,l_0)}{a,a_0}\Zpl{(l_0,r_0)}{a_0,b_0}\bl(\calA_1\br)
\Zpl{(r_0,r)}{b_0,b}
\\[1ex]
\le2M^2e^{-\veps(\neps^0-2)/f(2/\rho)}\Zpl{(l,r)}{a,b}
\end{gather*}
with sum running over $\rho\Hone\le a_0,b_0\le2T\Hone$.
Finally, taking $K_0\ge8$ such that
\begin{equation}\label{eq:K0large}
2M^2\le\exp\bl\{K_0\veps/4f(2/\rho)\br\}
\end{equation}
we immediately get
\[
\Zpl{(l,r)}{a,b}\bl(\calA_2\br)
\le\exp\Bl\{-\frac{\neps^0\veps}{2f(2/\rho)}\Br\}
\,\Zpl{(l,r)}{a,b}
\]
and thus \eqref{eq:OneDropBound}. It remains to verify 
\eqref{eq:RefinedFact1} and \eqref{eq:RefinedFact2}.

The proof of \eqref{eq:RefinedFact1} follows the argument of
Sect.~\ref{subsec:UpperBoundThm1.2}. Clearly, it is enough to show
that for some constant $c_1>0$ one has
\[
\Zpl{(l,r)}{a,b}\bl(\calA_2,a_0>2T\Hone\br)
\le\exp\bl\{-c_1T^{3/2}\br\}\Zpl{(l,r)}{a,b}\bl(\calA_2\br)\,.
\]
Using the definitions \eqref{eq:lrDef} and \eqref{eq:ArhoDef}, 
the partition function in the LHS of the previous display is bounded
above by
\[
\sum\Zpl{(l,l')}{a,a'}
\Zpl{(l',r')}{a',b'}\bl(\calA_{2\rho}(l',r'),a_0>2T\Hone\br)
\Zpl{(r',r)}{b',b}\,,
\]
where the sum runs over $0\le a',b'\le2\rho\Hone$ and $l'$, $r'$ such
that (recall \eqref{eq:lrDef})
\[
(l_0,r_0)\subset(l',r')\subseteq(l,r)\qquad\text{ \small and }\qquad
\Delta_{2\rho}=(l',r')\,. 
\]
Using Corollary~\ref{cor_HighDropletBound} the sum above can be
further majorated by
\[
e^{-c_1T^{3/2}}\sum\Zpl{(l,l')}{a,a'}
\Zpl{(l',r')}{a',b'}\bl(\calA_{2\rho}(l',r')\br)
\Zpl{(r',r)}{b',b}
\le e^{-c_1T^{3/2}}\Zpl{(l,r)}{a,b}\bl(\calA_2\br)\,.
\]
As a result, 
\[
\Zpl{(l,r)}{a,b}\bl(\calA_2,a_0\le2T\Hone,b_0\le2T\Hone\br)
\ge\bl(1-2e^{-c_1T^{3/2}}\br)\Zpl{(l,r)}{a,b}\bl(\calA_2\br)
\]
and it remains to choose $T$ large enough. The estimate
\eqref{eq:RefinedFact1} follows. 

Finally, we check \eqref{eq:RefinedFact2}. First, applying an obvious
extension of Lemma~\ref{lem_lowerboundPF}, we remove the field
$\lambda$ (as above, we put $b'=b-\rho\Hone$): for some
$M_1=M_1(\veps,\rho,T)$,  
\begin{equation}
\label{eq:F2M1Bound}
\frac{\Zpl{\Delta}{a,b}}{\Zpl{\Delta}{a,b'}}
\le M_1\frac{\Zpo{\Delta}{a,b}}{\Zpo{\Delta}{a,b'}}
= M_1\frac{\Ppo{\Delta}{a}\bl(X_\Delta=b\br)}
{\Ppo{\Delta}{a}\bl(X_\Delta=b'\br)}
\end{equation}
uniformly in $a$, $b$, $\Delta$ under consideration and all
$\lambda>0$ small enough. Here and below, $\Ppo{\Delta}{a}(\cdot)$
denotes the probability distribution of the $\Delta$-step random walk
starting from $a$ with transition probabilities $p(\cdot)$ restricted
to the set $\calI_{\Delta,+}$ of non-negative trajectories
(recall \eqref{eq:INp.Def}).

Now, denoting by $\ge0$ the wall constraint
$\BbbX\in\calI_{\Delta,+}$, we rewrite the last ratio as  
\begin{equation}
\label{eq:F2M2Bound}
\frac{\rmP_\Delta^a\bl(X_\Delta=b\mid\,\ge0\br)}
{\rmP_\Delta^a\bl(X_\Delta=b'\mid\,\ge0\br)}
=\frac{\rmP_\Delta^a\bl(X_\Delta=b,\ge0\br)}
{\rmP_\Delta^a\bl(X_\Delta=b',\ge0\br)}
\end{equation}
and observe that uniformly in $a$, $\Delta$ under consideration the
$\rmP_\Delta^a$-probability of the event $\ge0$ is uniformly positive. 
Thus, applying the standard argument (see, eg., \cite[\S11]{pB68},
\cite[\S9]{pB99}) one
deduces that, uniformly in $a$, $b$, $\Delta$, and $\lambda$ under
consideration, the ratio in \eqref{eq:F2M2Bound} is bounded above by a
positive constant $M_2=M_2(\veps,\rho,\Delta,T)$.  
The estimate\eqref{eq:RefinedFact2} follows from \eqref{eq:F2M1Bound}
and  \eqref{eq:F2M2Bound}.

The proof of the lemma is finished.
\end{proof}

Next, we fix $\rho$, $K$, and $\veps$ as in the proof above, use the
$K$-blocks decomposition and introduce the labels
(cf.~\eqref{eq:YLabels}) 
\begin{equation}
\label{eq:ULabels}
U_m=U_m^\rho(X)
=\ind{\{\max_{j\in\calI_m}X_j<4\rho\Hone \}}\,.
\end{equation}
Denote $n_U=\left|\bl\{m:U_m=1\br\}\right|$.

\begin{lemma}\label{lem:nUBound}\sl
For any $\rho>0$ there exist positive constants $\lambda_0$, $\gamma_0$,
$K_0$, $c$, and $C$ 
such that for any $\lambda\in(0,\lambda_0)$, $\gamma\in(0,\gamma_0)$, 
$K\ge K_0$ and all $N\ge3K\Hone(\lambda)^2$ we have
\begin{equation}
\label{eq:nUBound}
\PNpl{a,b}\bl(n_U\le\gamma n_K\br)
\le C\exp\Bl\{-cN \Hone(\lambda)^{-2}\Br\}
\end{equation}
uniformly in $a$, $b\in\bl(0,2\rho\Hone(\lambda)\br)$.
\end{lemma}

\begin{proof}
Our argument is similar to that of Sect~\ref{sec:LowerBound}. First,
taking $K=3K_0$ in Lemma~\ref{lem:nYBound}, we split all $K$-blocks
into triples of consecutive blocks (neglecting the non-complete last
triple if there is one) and call an index $m$ regular if 
\[
Y_{3m+1}=Y_{3m+3}=0\,.
\]
Using the previous Lemma with $\alpha=1/9$, we deduce that with
probability not smaller than 
$1-\exp\bl\{-N/(72f(2/\rho)\Hone(\lambda)^{2})\br\}$ there
are at least $2n_K/9$ regular indices $m$.

For each such $m$ we define
\begin{gather*}
l=\max\Bl\{j\in\calI_{3m+1}:X_j<2\rho\Hone(\lambda) \Br\},
\quad\quad X_l=a,\\
r=\min\Bl\{j\in\calI_{3m+3}:X_j<2\rho\Hone(\lambda) \Br\},
\quad\quad X_r=b.
\end{gather*}
Now, $K\veps\Hone(\lambda)^2=3K_0\veps\Hone(\lambda)^2
\le r-l\le3K\veps\Hone(\lambda)^2$
and thus, using Corollary~\ref{cor:RemovingLambda} we get
\begin{equation}\label{eq:pK.def}
\begin{split}
\Ppl{(l,r)}{a,b}\bl(U_{3m+2}=1\br)
&\ge ce^{-C (r-l)/(2\rho\Hone(\lambda))^{2}}
\\[1ex]
&\ge ce^{-3CK\veps/(4\rho^2)}=:p_{\rho,K\veps}>0\,.
\end{split}
\end{equation}
Therefore, on average there are at least $2p_{\rho,K\veps}n_K/9$
indices $m$ whose labels satisfy
$U_{3m+2}=1$. By a standard large deviation bound we get
\[
\PNpl{}\Bl(n_U<\frac{p_{\rho,K\veps}}{9}n_K\Br)\le Ce^{-cn_K}
\]
and the lemma is proved.
\end{proof}

\subsection{Coupling}\label{sec:coupling}

With $\rho>0$ fixed as above and $\lambda>0$ denote (recall
\eqref{eq:Hgam.Def}) 
\begin{equation}\label{eq:H.def}
H=H(\rho,\lambda):=4\rho\Hone(\lambda).
\end{equation}
Let $X$, $Y$ be two independent trajectories of our process and let
$\PPNpl{\abxy}(\ppunkt)$ denote their joint distribution
(recall \eqref{eq:INp.Def}--\eqref{eq:PNpl.Def}),
\[
\PPNpl{\abxy}(\ppunkt)=\PNpl{a_x,b_x}(\punkt)
\otimes\PNpl{a_y,b_y}(\punkt)
\]
with the shorthand notation $\PPNpl\nul(\punkt)$ if
$a_x=b_x=a_y=b_y=0$. Everywhere in this section we shall consider
only boundary conditions satisfying $0\le\abxy\le H$.
For a set of indices $A\subseteq[0,N]\cap\BbbZ$, let
\[
\calN_A=\calN_A(X,Y)=\Bl\{\every j\in A, X_j\neq Y_j\Br\}
\]
be the event ``trajectories $X$ and $Y$ do not intersect within
the set $A$''.
Our main observation is that, with probability going to one, any two
independent trajectories of our RW meet within a time interval of
order at most $H^2=O\bl(\Hone(\lambda)^2\br)$: 

\begin{lemma}\sl
There exist positive constants $\lambda_0$, $C$, $c$, and $\rho_0$
such that the inequality
\begin{equation}\label{eq:nocrossN}
\PPNpl\nul\Bl(\NN\Br)\le Ce^{-cN/H^2}
\end{equation}
holds uniformly in $0<\lambda\le\lambda_0$, $0<\rho\le\rho_0$ and
$N\ge H^2$.
\end{lemma}

\begin{proof}
Consider the decomposition into $K$-blocks $\calI_m$ of length
$K\veps\Hone(\lambda)^2$ described in the previous section and denote
\[
n_U^2=\left|\Bl\{m:U_m(X)=U_m(Y)=1\Br\}\right|\,
\]
with labels $U_m(\punkt)$ defined as in \eqref{eq:ULabels}. Following
the  proof of Lemma~\ref{lem:nUBound} with $\alpha=1/9$ we deduce that
with probability not smaller than 
\[
1-2\exp\Bl\{-\frac N{72f(2/\rho)\Hone(\lambda)^{2}}\Br\}
\]
there are on average at least $\bl(1/3-2\alpha\br)n_K=n_K/9$ regular
indices $m$ that are common for both $X$ and $Y$. As a result (recall
\eqref{eq:pK.def}), 
\begin{equation}\label{eq:n2U.bound}
\PNpl{}\Bl(n_U^2<\frac{(p_{\rho,K})^2}{18}n_K\Br)\le C_1e^{-c_1n_K}
\le C_1e^{-c_2N/H^2}\,,
\end{equation}
that is, with high probability there is a positive fraction of blocks
$\calI_m$ for which the event 
\begin{equation}\label{eq:Dm.def}
D_m=\Bl\{U_{m}(X)=U_{m}(Y)=1\Br\}
=\Bl\{\every j\in\calI_m,0\le X_j,Y_j<4\rho\Hone \Br\}\,
\end{equation}
is realized. By taking each second such $K$-block we construct a
disjoint collection $\calK$ of $K$-blocks possessing property
\eqref{eq:Dm.def}. The collection $\calK$ has the following 
important properties to be used in the sequel:
\begin{quote}
{\sl
1) with probability at least $1-C_1\exp\bl\{-c_2NH^{-2}\br\}$, there are
   no less than $p_{\rho,K}^2n_K/36$ blocks in $\calK$;\\
2) conditioned on $\bl\{U_m\br\}$ and on the configuration in the
   complement of $K$-blocks from $\calK$, distributions inside
   individual $K$-blocks are independent.
}
\end{quote}

Another important ingredient of our argument is the following
observation:

\begin{lemma}\label{lem:nocross}\sl
Let $H=H(\rho,\lambda)$ be as in \eqref{eq:H.def} and let
$D=D_N(X)\cap D_N(Y)$, where 
\begin{equation}\label{eq:DN.def}
D_N(Z)=\Bl\{0\le Z_j\le H, \every j=1,\ldots,N-1\Br\}.
\end{equation}
Then there exist positive constants $\lambda_0$ and $c_3$ such that the
estimate 
\begin{equation}\label{eq:nocross.bound}
\max_{\abxy}\PPNpo{\abxy}\Bl(\NN\mid D\Br)\le e^{-c_3NH^{-2}},
\end{equation}
holds uniformly in $0<\lambda\le\lambda_0$, in $N\ge H^2$ and in
boundary conditions $0\le \abxy\le H$. 
\end{lemma}

We postpone the proof of Lemma~\ref{lem:nocross} till the end of this
section and deduce our main estimate \eqref{eq:nocrossN} first.
Combining Corollary~\ref{cor:RemovingLambda} with the inequality
\eqref{eq:nocross.bound} and using the estimate 
\[
\lambda V(4\rho\Hone)N\le
\frac{2\rho N}{\Hone^2}\lambda V(2\Hone)\Hone^2
=\frac{32\rho^3N}{H^2}
\]
we obtain the uniform bound (similarly as in~\eqref{eq:RemovingLambda}) 
\begin{equation}\label{eq:nocross.1}
\max_{\abxy}\PPNpl\abxy\Bl(\NN\mid D\Br)
\le e^{-(c_3-64\rho^3)NH^{-2}}\le e^{-c_3N/2H^2}
\end{equation}
provided only $0<\rho\le\rho_0$ with $128\rho_0^2<c_3$. Now, using the
bound 
\[
\begin{split}
\PPNpl\nul\Bl(\NN\Br)\le\PPNpl\nul\Bl(n_U^2&<\frac{(p_{\rho,K})^2}{18}n_K\Br)\\
&+\PPNpl\nul\Bl(\NN\bigm||\calK|\ge\frac{(p_{\rho,K})^2}{36}n_K\Br),
\end{split}
\]
``freezing'' the joint configuration $(X,Y)_{\calK^c}$ in all blocks
that do not belong to the collection $\calK$ and 
using the estimate \eqref{eq:nocross.1} for all blocks from $\calK$,
we bound the last term by
\begin{equation}\label{eq:nocross.2}
\begin{split}
\max_{(X,Y)_{\calK^c}}\PPNpl\nul&\Bl(\NN\bigm|(X,Y)_{\calK^c},
|\calK|\ge\frac{(p_{\rho,K})^2}{36}n_K\Br)\\
&\le\exp\Bl\{-\frac{c_3|\calI_m|}{2H^2}\cdot\frac{(p_{\rho,K})^2}{36}n_K\Br\}
\le e^{-c_4NH^{-2}}.
\end{split}
\end{equation}
Averaging this inequality over $(X,Y)_{\calK^c}$ and combining the
result with \eqref{eq:n2U.bound}, we obtain the target estimate
\eqref{eq:nocrossN}. 
\end{proof}

\begin{proof}[Proof of Lemma~\ref{lem:nocross}]
We use again a blocking argument. Fix two positive constants $\veps$
and $\Delta$ such that $\veps\le1$, $\Delta\le1/8$ and assume
w.l.o.g.\ that $\veps H^2$ is integer. Split the interval $[0,N]$ into
blocks $\calI_m$ of length $4\veps H^2$, 
$m=1,2,\dots,\lfloor N/(4\veps H^2)\rfloor=:\neps$, and for each such
block introduce the ``crossing event''
\[
\Am(X,Y)=\Amup(X)\cap\Amdn(Y),
\]
where (see Fig.~\ref{fig:crossing})
\begin{equation}\label{eq:Aevent}
\begin{gathered}
\Amup(Z)=D_m(Z)\cap \hatDom(Z)\cap\Bmup(Z),\\[1ex]
\Amdn(Z)=D_m(Z)\cap \hatDom(Z)\cap\Bmdn(Z)
\end{gathered}
\end{equation}
with $D_m=D_{\calI_m}$ defined as in \eqref{eq:Dm.def},
$\hatDom=\hatD_{\calI_m^o}$ given by
\begin{gather*}
\hatD_J(Z)=\Bl\{\frac14H\le Z_j\le\frac34H,\every j\in J\Br\},
\qquad\every J\subseteq[0,N],\\ 
\calI_m^o=\Bl\{j\in\calI_m:(4m-3)\veps H^2<j<(4m-1)\veps H^2\Br\}
\end{gather*}
and, finally, 
\begin{gather*}
\Bmup(Z)=\Bl\{Z_{(4m-3)\veps H^2}\in\CmD,Z_{(4m-1)\veps H^2}\in\CpD\Br\},
\\[1ex]
\Bmdn(Z)=\Bl\{Z_{(4m-3)\veps H^2}\in\CpD,Z_{(4m-1)\veps H^2}\in\CmD\Br\}
\end{gather*}
with
\[
\CpD=\Bl[\frac12H,\frac{1+\Delta}2H\Br],\qquad
\CmD=\Bl[\frac{1-\Delta}2H,\frac12H\Br].
\]

\begin{figure}[ht]
\bigskip
\centerline{\includegraphics[height=35mm]{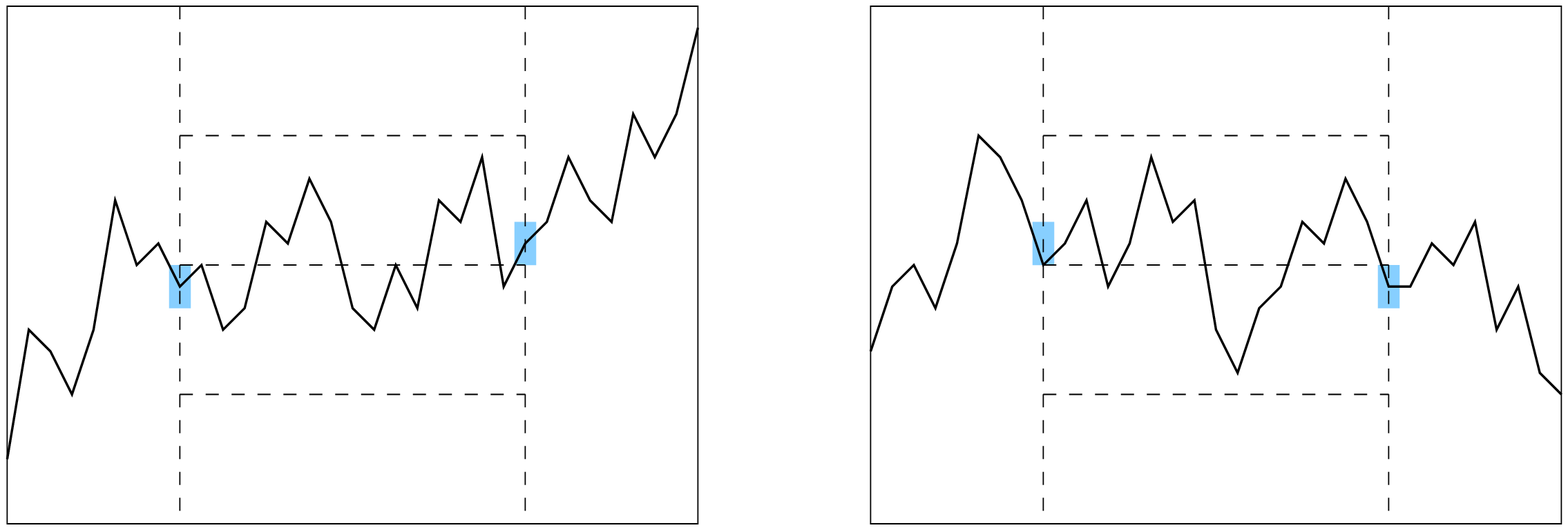}}
\bigskip
\caption{Crossing events $\Amup$ (left) and $\Amdn$ (right)}
\label{fig:crossing}
\end{figure}

Our argument below is based upon the following three facts:
{\sl
\begin{quote}
F1) conditioned on the ``boundary'' values
\[
\Bl\{(X_j,Y_j),j=4m\veps H^2, m=1,2,\dots,\neps\Br\},
\] 
the trajectories $X$ and $Y$ behave independently in different blocks; 

F2) for all $\calI_m$, uniformly in ``boundary'' values,  the crossing
   event $\Am(X,Y)$ occurs with positive probability:
\end{quote}
\begin{equation}\label{eq:pone.bound}
\min_{\abxy}\PPmpo{\abxy}\Bl(\Am(X,Y)\mid D_m\,\Br)\ge p_1>0
\end{equation}
\begin{quote}
F3) conditioned on $\Am(X,Y)$, the trajectories $X$ and $Y$ intersect
   in $\calI_m^o$ with positive probability:
\end{quote}
\begin{equation}\label{eq:ptwo.bound}
\min_{\abxy}\PPmpo{\abxy}\Bl(\noNm\mid\Am(X,Y)\Br)\ge p_2>0
\end{equation}
\begin{quote}
where $\noNm$ denotes the intersection event in the central part of
the block $\calI_m$:
\end{quote}
\[
\noNm=\Bl\{\some j\in\calI_m^o:
X_j=Y_j\in\Bl[\frac18H,\frac78H\Br]\Br\}.
\]
}

Indeed, let $\calA$ be the collection of indices $m$ for which the
event $\Am(X,Y)$ occurs. Denoting by $n_{\sf cross}$ the cardinality of
$|\calA|$, we use the Markov property and the standard large deviation
bound to get
\begin{equation}\label{eq:ncross.bound}
\max_{\abxy}\PPNpo\abxy\Bl(n_{\sf cross}<\frac{p_1}2\neps\mid D_m\Br)
\le e^{-c_5\neps}
\end{equation}
with some $c_5>0$; thus, it remains to consider the case 
$n_{\sf cross}=|\calA|\ge p_1\neps/2$. Now, the bound
\eqref{eq:ptwo.bound} and the conditional independence of blocks
$\calI_m$ imply that
\[
\begin{split}
\PPNpo\abxy\Bl(\NN\mid\calA,D_m\Br)
&\le\prod_{m\in\calA}\PPmpo\abxym\Bl(\Nm\mid\Am(X,Y),D_m\Br)\\
&\le(1-p_2)^{|\calA|}
\end{split}
\]
uniformly in all such $\calA$ and boundary conditions $\abxym$ thus
giving the uniform upper bound
\begin{equation}\label{eq:ncrossAD.bound}
(1-p_2)^{p_1\neps/2}\le e^{-c_6\neps}\,.
\end{equation}
Finally, \eqref{eq:nocross.bound} follows immediately from
\eqref{eq:ncross.bound} and \eqref{eq:ncrossAD.bound}.

Thus, it remains to check the properties F1)--F3).
The Markov property F1) being obvious, we need only to prove the
inequalities \eqref{eq:pone.bound} and \eqref{eq:ptwo.bound}. 

To check \eqref{eq:pone.bound} we proceed as in the proof of
Lemma~\ref{lem_lowerboundPF}. By conditional independence, it is
enough to show that for some $p_3>0$ one has
\begin{gather*}
\min_{\abx}\Pmpo\abx\bl(\Amup(X)\mid D_m\br)\ge p_3,\\
\min_{\aby}\Pmpo\aby\bl(\Amdn(Y)\mid D_m\br)\ge p_3.
\end{gather*}
We shall verify the first of these inequalities, the second follows
from analogous considerations.

First, for any boundary conditions $a,b$ we rewrite (recall
\eqref{eq:Aevent})
\begin{align*}
\Pmpo{a,b}\bl(\Amup\mid D_m\br)
&=\Pmpo{a,b}\bl(\hatDom\cap\Bmup\mid D_m)\\
&=\Pmpo{a,b}\bl(\hatDom\mid\Bmup\cap D_m)\cdot
\Pmpo{a,b}\bl(\Bmup\mid D_m).
\end{align*}
Next, using the Markov property and the usual invariance principle (for
unconditioned RWs), we get
\[
\begin{split}
\Pmpo{a,b}\bl(\hatDom\mid\Bmup\cap D_m)
\ge&\min_{x\in\CmD,y\in\CpD}
\Bl[\mathrm{P}_{2\veps H^2}^{x,y}\bl(\hatDom\br)-
    \mathrm{P}_{2\veps H^2}^{x,y}\bl(\overline{D}_{\calI^o_m}\br)\Br]\\
\ge\min_{x\in\CmD,y\in\CpD}
&\mathrm{P}_{2\veps H^2}^{x,y}\Bl(\max|X_j-x|\le(1-2\Delta)H/4\Br)
\\-\max_{x\in\CmD,y\in\CpD}
&\mathrm{P}_{2\veps H^2}^{x,y}\Bl(\max|X_j-x|>(1-\Delta)H/2\Br)
\ge p_4
\end{split}
\]
where, $\overline{D}_{\calI^0_m}$ denotes the complement of the event
$D_{\calI^o_m}$ (cf.~\eqref{eq:Dm.def}),
\[
D_{\calI^o_m}=\Bl\{\every j\in\calI^o_m,0\le X_j,Y_j<4\rho\Hone \Br\}\,,
\]
and for fixed $\Delta\in(0,1/8]$, the bound $p_4=p_4(\veps,\Delta)$ is
positive uniformly in all $\lambda$ and $\veps$ small enough. 
On the other hand, a slight modification of the proof of
Lemma~\ref{lem_lowerboundPF} implies the inequality
\[
\Pmpo{a,b}\bl(\Bmup\mid D_m)\ge p_5
\]
with some $p_5=p_5(\veps,\Delta)>0$ uniformly in all $\lambda$ small
enough. The estimate \eqref{eq:pone.bound} follows.

Finally, we verify the inequality \eqref{eq:ptwo.bound}. Morally, our
argument is based upon the following observation: on the event
$\Am(X,Y)$ there is $j\in\calI_m^o$ satisfying
\begin{equation}
X_j-Y_j\le0\le X_{j+1}-Y_{j+1};
\end{equation}
since $X$ and $Y$ are independent processes whose jumps have the same
distribution of finite variance $\sigma^2>0$, the bound
\begin{equation}\label{eq:bounded.jump}
\min\Bl\{|X_j-Y_j|,|X_{j+1}-Y_{j+1}|\Br\}\le R\sigma
\end{equation}
holds with positive probability provided only the absolute constant
$R>0$ is chosen large enough; finally, thanks to the aperiodicity
property \eqref{eq:aperiodic}, conditioned on the event
\eqref{eq:bounded.jump}, the trajectories $X$ and $Y$ meet with
positive probability within $A\lceil R\sigma\rceil$ steps.

We sketch the main steps of the argument.
In view of the Markov property, it is sufficient to show that 
\[
\min_{x\in\CmD,y\in\CpD}
\mathrm{P}_{\calI_m^o}^{\,x,y}
\Bl(\some j\in\calI_m^o,\,|X_j-Y_j|\le R\sigma\mid \hatDom\Br)
\]
is uniformly positive for all $\lambda$ small enough. Clearly, the
minimum above is bounded below by the expression
\[
{\sf Pr}\Bl(\bl|\xi-\eta\br|\le R\sigma\Br)
-2\max\mathrm{P}_{\calI_m^o}^{\,x,y}
\Bl(\max_{j\in\calI_m^o}|X_j-x|>\frac{1-2\Delta}4H\Br),
\]
where $\xi$ and $\eta$ are i.i.d.~r.v.\ with the basic distribution
$p(\punkt)$ and $\max$ is taken over all $x\in\CmD$, $y\in\CpD$. 
As $R\to\infty$, the first term approaches $1$, whereas the second
vanishes asymptotically as $\veps\to0$, uniformly in
$\Delta\in(0,1/8]$ and in all $\lambda>0$ small enough. Thus, for some
$p_6=p_6(\veps,\Delta,R)>0$ we get
\[
\min_{\abxym}\PPmpo{\abxym}\Bl(\some j\in\calI_m^o,\,|X_j-Y_j|\le R\sigma
\mid \hatDom\Br)\ge p_6
\]
and thanks to the aperiodicity property \eqref{eq:aperiodic} the
trajectories $X$ and $Y$ have a positive probability to meet within the
time interval $J_0=[j_0,j_0+A\lceil R\sigma\rceil]$:
\[
\PPmpo{\abxym}
\Bl(\some j\in J_0,\, X_j=Y_j\in[0,H]
\Bigm|
\Am(X,Y),|X_{j_0}-Y_{j_0}|\le R\sigma\Br)\ge p_7
\]
with some $p_7>0$, uniformly in boundary conditions
$0\le\abxym\le H$ and in
positive $\lambda$ small enough. This implies the
estimate \eqref{eq:ptwo.bound}. 

The proof of the lemma is complete.
\end{proof}


\subsection{Relaxation to equilibrium}

It is an immediate corollary of \eqref{eq:nocrossN}: just consider the
initial RW and another one started at equilibrium. By the coupling
inequality \cite[pg.~12]{tL92} 
the total variance distance between the distribution of our
RW after $N$ steps and the equilibrium measure is bounded above by
the LHS of \eqref{eq:nocrossN}.

\subsection{Inverse correlation length}

Let $X$ be our RW and $Y$ its independent copy; we have:
\begin{equation}
\label{eq:covariance}
\Cov(X_i,X_j)=\frac12\epl\bl[(X_i-Y_i)(X_j-Y_j)\br]\,,
\end{equation}
where $\epl$ is the expectation w.r.t. the limiting measure and $\Cov$
is the corresponding covariance; denote by $A$ the event that both
RW's $X$ and $Y$ intersect between $i$ and $j$.
According to the above, the probability of the complement
$\bar A$ of $A$ is bounded above by the RHS of \eqref{eq:nocrossN}: 
\[
\ppl(\bar A)\le C\exp\bl\{-c\,|i-j|\,\Hone^{-2}\br\}.
\]
Moreover, by symmetry of the RHS of \eqref{eq:covariance} on the event
$A$, we have \[
\epl\bl[(X_i-Y_i)(X_j-Y_j)\ind A\br]=0.
\]
Consequently, for any $p>1$,
\[
\begin{split}
2\Cov(X_i,X_j)&=\epl\bl[(X_i-Y_i)(X_j-Y_j)\ind{\bar A}\br]\\[1ex]
&\le\epl\bl[(X_iX_j+Y_iY_j)\ind{\bar A}\br]\\[1ex]
&\le 2\bl(\epl\bl[(X_iX_j)^{p}\br]\br)^{1/p}\bl(\ppl(\bar A)\br)^{(p-1)/p}.
\end{split}
\]
However, for any $p$, $1<p<21/8$ we get (recall
Corollary~\ref{cor_HighMomentUpperBound})
\[
\epl{}\bl(\,X_{N_2}\,\br)^{2p}\le C(p)\,\Hone^{2p+1}\,,
\]
and thus, by the Cauchy-Schwarz inequality, 
\[
\epl\bl[(X_iX_j)^p\br]\le\Bl[\epl\bl((X_i)^{2p}\br)
\epl\bl((X_j)^{2p}\br)\Br]^{1/2}
\le C\Hone^{2p+1}
\]
leading to 
\[
\Cov(X_i,X_j)\le C\Hone(\lambda)^{2+1/p}
\exp\bl\{-c\,|i-j|\,\Hone^{-2}\br\}.
\]
Finally, take $p=2$.
\qed

\appendix
\section{Small droplet bound}\label{sec:SmallDropletBound}
Our aim here is to prove the small droplet
bound---Lemma~\ref{lem:SmallDropletBound}. The key step of our
argument will be based upon the following, having an independent
interest, conditional Chebyshev inequality for maximum.

\begin{lemma}\label{lem:ConditionalChebyshevForMaximum}\sl
Let $S_0=0$, $S_k=\xi_1+\dots+\xi_k$, $k\ge1$, be the random walk
generated by a sequence $\xi_1$, $\xi_2$, \dots of i.i.d. random
variables such that $\bfE\xi=0$, $\bfE\xi^2=\sig^2<\infty$. 
Let $D>0$ be an arbitrary constant and, for any $m\ge1$, let $d_m$
satisfy $\bfP(S_m=d_m)>0$ and  $|d_m|\le D$.
Then there exists a positive constant $c=c(D)$ such that the inequality
\begin{equation}\label{eq:ConditionalChebyshevForMaximum}
\bfP\bl(\max_{0<k<m}S_k>M\mid S_m=d_m\br)\le c\,\frac{m^{3/2}}{M^2}
\end{equation}
holds for all $m\ge2$.
\end{lemma}

\begin{proof}
Since $D>0$ is a finite constant, the local limit theorem \cite{bdG62}
implies that for some $c_1=c_1(D)>0$
\begin{equation}\label{eq:LLTLowerBound}
\bfP\bl(S_m=d_m\br)\ge
\frac{c_1}{\sqrt{2\pi\sig^2m}}\exp\Bl\{-\frac{(d_m)^2}{2\sig^2m}\Br\}
\end{equation}
uniformly in $m\ge1$ and $|d_m|\le D$.

On the other hand, by the Etemadi (see, eg, \cite[pg.~256]{pB99}) and
Chebyshev inequalities,
\[
\bfP\bl(\max_{0<k<m}S_k>M\br)\le3\max_{0<k<m}\bfP\bl(S_k>M/3\br)
\le3\max_{0<k<m}\frac{\sig^2k}{(M/3)^2}=\frac{27\sig^2m}{M^2}\,.
\]
The target bound \eqref{eq:ConditionalChebyshevForMaximum} follows
immediately from the last two displays and the assumption $|d_m|\le D$.
\end{proof}

\begin{proof}[Proof of Lemma~\ref{lem:SmallDropletBound}]
Let first $m\le M^{7/6}$. Then, the Conditional Chebyshev inequality
\eqref{eq:ConditionalChebyshevForMaximum} gives
\[
\bfP\bl(\max_{0<k<m}S_k>M\mid S_m=d_m\br)\le C_1\,\frac{m^{3/2}}{M^2}
\le \frac{C_1}{M^{1/4}}\,.
\]
Let now $m$ satisfy $M^{7/6}\le m\le\zeta M^2$ (and thus
$m\to\infty$). Since $D>0$ is finite, 
it follows from the main result in \cite{wdK76} that 
\[
\bfP\bl(\max_{0<k<m}S_k>M\mid S_m=d_m\br)\le 
C_2\exp\Bl\{-C_3\frac{M^2}m\Br\}\le C_2\exp\Bl\{-\frac{C_3}\zeta\Br\}
\]
if only $M$ is large enough, $M\ge M_0$. The small droplet bound
\eqref{eq:SmallDropletBound} follows, provided $\zeta>0$ is chosen
small enough.
\end{proof}

Next, we present a simple one-point analogue of
Lemma~\ref{lem:ConditionalChebyshevForMaximum}. 

\begin{lemma}\label{lem_OnePointConditionalChebyshev}\sl
Under the conditions of
Lemma~\ref{lem:ConditionalChebyshevForMaximum}, there is a positive
constant $\bar c$ depending on $D$ and the distribution of $\xi$
only, such that 
\[
\max_{0<k<m}\bfP\bl(S_k>M+D\mid S_m=d_m\br)
\le\bar c\,\frac {m^{5/2}}{M^4}
\]
for all $m\ge2$.
\end{lemma}

\begin{proof}
Using the independence of jumps and the Chebyshev inequality, we get
\begin{equation}\label{eq:OnePointChebyshev}
\begin{split}
\bfP\bl(S_k>M+D,S_m=d_m\br)&\le\bfP\bl(S_k>M\br)\bfP\bl(S_m-S_k<-M\br)
\\[1ex]
&\le\frac{k(m-k)}{M^4}\sig^4\le\frac{m^2}{4M^4}\sig^4\,.
\end{split}
\end{equation}
Combining this estimate with the lower bound~\eqref{eq:LLTLowerBound},
we deduce the result.
\end{proof}

Finally, we present a stronger version of the previous claim.

\begin{lemma}\label{lem_ConditionalChebyshev}\sl
Under the conditions of
Lemma~\ref{lem:ConditionalChebyshevForMaximum}, there is a positive
constant $\tilde c$ depending on $D$ and the distribution of $\xi$
only, such that 
\[
\max_{0<k<m}\bfP\bl(S_k>M+D\mid S_m=d_m\br)\le \tilde c\,\frac m{M^2}\,.
\]
\end{lemma}

\begin{proof}
We start by observing that if $\xi_1$, $\xi_2$, \dots, are i.i.d. random
variables and $d_m$ is chosen such that $\bfP(S_m=d_m)>0$, then the
variables $\eta_j$ defined via $\eta_j=(\xi_j\mid S_m=d_m)$ are
exchangeable. As a result \cite[\S24]{pB68}, for any $k$, 
$1\le k\le m$, 
\begin{equation}\label{eq:ExchangeableRelations}
\begin{split}
\bfE\bl(\,S_k\mid S_m=d_m\,\br)&
=k\,\bfE\bl(\,\xi_1\mid S_m=d_m\,\br)=\frac{ka}n\,,
\\
{\mathbf {Var}}\bl(\,S_k\mid S_m=d_m\,\br)&=\frac{k(m-k)}{m-1}\,
{\mathbf {Var}}\bl(\,\xi_1\mid S_m=d_m\,\br)\,.
\end{split}
\end{equation}

Our next observation formalizes an intuitively obvious fact that for
large $m$ the variable $\xi_k$ becomes asymptotically independent of
$S_m$ and thus the variances of $\eta_1$ and $\xi_1$ are close to each
other. We shall restrict ourselves to the case of integer-valued
variables $\xi$ having zero mean and the variance $\bfE\,\xi^2=\sig^2$.
\begin{quote}\sl
For any finite $D>0$ there exists $m_0$, depending only on $D$ and the
distribution of $\xi$, such that the inequality
\end{quote}
\begin{equation}\label{eq:ConditionalVarianceBound}
\bfE\bl(\xi_1{}^2\mid S_m=d\,\br)\le4\bfE\,\xi_1{}^2=4\sig^2
\end{equation}
\begin{quote}\sl
holds uniformly in $m\ge m_0$ and $|d|\le D$.
\end{quote}
To check \eqref{eq:ConditionalVarianceBound}, we observe that the
characteristic function of $\eta=(\xi\mid S_m=d)$ equals
\[
\bfE\bl(e^{is\xi}\mid S_m=d\,\br)=
\frac{\int\phi(t+s)\phi^{m-1}(t)e^{-itd}\,dt}
{\int\phi^m(t)e^{-itd}\,dt}\,,
\]
where $\phi(t)$ is the unconditional characteristic function of $\xi$,
$\phi(t)=\bfE e^{it\xi}$, and the integration goes over an interval of
periodicity of $\phi(t)$. Consequently,
\[
\bfE\bl(\,\xi^2\mid S_m=d\,\br)=
\frac{-\int\phi''(t)\phi^{m-1}(t)e^{-itd}\,dt}
{\int\phi^m(t)e^{-itd}\,dt}\,.
\]

According to \eqref{eq:LLTLowerBound}, we have
\[
\int\phi^m(t)e^{-itd}\,dt=\bfP(S_m=d)\ge
\frac{1}{2\sqrt{2\pi\sig^2m}}\exp\Bl\{-\frac{d^2}{2\sig^2m}\Br\}
\]
uniformly in $|d|\le D$ and all $m\ge m_1$ with $m_1$ large enough.
Analogously, applying the standard Laplace method to the integral in
the numerator (see, eg, \cite{ngDB58,W41}), we get
\[
\Bl|\int\phi''(t)\phi^{m-1}(t)e^{-itd}\,dt\Br|\le
\frac{2}{\sqrt{2\pi\sig^2m}}\exp\Bl\{-\frac{d^2}{2\sig^2m}\Br\}\,
\bfE\,\xi^2\,,
\]
uniformly in $|d|\le D$ and all $m\ge m_2$ with $m_2$ large enough.
The bound \eqref{eq:ConditionalVarianceBound} follows from the last
two displays.

Next, we combine \eqref{eq:ExchangeableRelations} and
\eqref{eq:ConditionalVarianceBound} to deduce that, uniformly in
$|d|\le D$ and all $m\ge m_0$ with $m_0$ large enough, the inequality
\[
\bfE\bl(\,S_k{}^2\mid S_m=d\,\br)\le\Bl(\frac{kd}m\Br)^2
+\frac{k(m-k)}{m-1}\,4\bfE\,\xi{}^2\le d^2+\frac{m^2}{m-1}\bfE\,\xi{}^2
\]
holds for all $k\in\bl\{1,\dots,m\br\}$. By Chebyshev,
\[
\bfP\bl(S_k>M+D\mid S_m=d_m\br)\le \frac{C_1m}{M^2}
\]
for all such $m$.

It remains to consider $m\le m_0$. Denoting 
\begin{equation}\label{eq:FiniteDropletProbability}
p(m_0,D)=\min_{m\le m_0,|d|\le D}
\Bl\{\bfP(\,S_m=d\,): \bfP(\,S_m=d\,)>0\Br\}>0\,,
\end{equation}
we immediately get, via Chebyshev,
\[
\bfP\bl(S_k>M+D\mid S_m=d_m\br)\le\frac{\bfP(S_k>M+D)}{p(m_0,D)}
\le\frac{C_2m}{M^2}\,.
\]
The proof is finished.
\end{proof}

\end{document}
\bibitem{rD}
Durrett, R.: 
paper in Ann. Probab.

\bibitem{tL} Liggett paper

\end{thebibliography}
\end{document}


\end{document}